\documentclass[10pt,onecolumn]{elsart3p}

\usepackage{graphicx,amssymb,epsfig,amsmath}
\usepackage[square,comma,sort,numbers]{natbib}

\numberwithin{equation}{section}

\newcommand{\operL}{{\mathcal L}}
\newcommand{\seB}{{\mathbb B}}
\newcommand{\er}{\mathbb R}
\newcommand{\caS}{{\mathcal S}}
\newcommand{\caR}{{\mathcal R}}
\newcommand{\D}{\textrm{d}}
\newcommand{\boldN}{\mathbf{N}}

\journal{Journal of Computational Physics}
\begin{document}

\begin{frontmatter}

\title{Projective and Coarse Projective Integration for Problems
with Continuous Symmetries}

\author[label1]{M.E. Kavousanakis \corauthref{cor1}},
\ead{mihkavus@chemeng.ntua.gr}
\author[label2]{R. Erban},
\ead{erban@maths.ox.ac.uk}
\author[label1]{A.G. Boudouvis},
\ead{boudouvi@chemeng.ntua.gr}
\author[label3,label4]{C.W. Gear},
\ead{wgear@princeton.edu}
\author[label3,label5]{I.G. Kevrekidis}
\ead{yannis@princeton.edu}

\address[label1]{National Technical University of Athens,
School of Chemical Engineering, 9 Heroon Polytechniou St.,
Zographos, Athens, GR-15780, Greece}
\address[label2]{University of Oxford, Mathematical Institute, 24-29 St. Giles',
Oxford, OX1 3LB, United Kingdom}
\address[label3]{Princeton University,
Department Of Chemical Engineering,
Engineering Quadrangle, Olden Street, Princeton, NJ 08544, USA}
\address[label4]{NEC Research Institute,retired}
\address[label5]{Princeton University, Program in Applied and Computational Mathematics,
Princeton, NJ 08544, USA}

\corauth[cor1]{Address: School of Chemical Engineering, 9 Heroon Polytechniou St., Zographos, Athens, Greece GR-15780,
               Tel.: +30 210 7723290, Fax: +30 210 7723155}

\begin{abstract}
Temporal integration of equations possessing
continuous symmetries
(e.g. systems with translational invariance associated with traveling solutions
and scale invariance associated with self-similar solutions)
in a ``co-evolving'' frame (i.e. a frame which is co-traveling, co-collapsing or
co-exploding with the evolving solution) leads to improved accuracy
because of the smaller time derivative in the new spatial frame.
The slower time behavior permits the use of {\it projective} and {\it coarse projective} 
integration with longer projective steps in the computation of the time evolution 
of partial differential equations and multiscale systems, respectively.
These methods are also demonstrated to be effective for systems which only 
approximately or asymptotically possess continuous symmetries.
The ideas of projective integration in a co-evolving frame are illustrated 
on the one-dimensional, translationally invariant Nagumo partial differential 
equation (PDE). A corresponding kinetic Monte Carlo model, motivated from the
Nagumo kinetics, is used to illustrate the coarse-grained method.
A simple, one-dimensional diffusion problem is used to illustrate the scale invariant 
case. 
The efficiency of projective integration in the 
co-evolving frame for both the macroscopic diffusion PDE and 
for a random-walker particle based model is again demonstrated.
\end{abstract}


\end{frontmatter}

\section{Introduction}
Projective and coarse projective integration have been recently proposed as effective
methods for the computation of long time behavior in complex multiscale
problems \cite{Gear:2001:PIM,Gear:2003:TPS,Gear:2003:PMS,Kevrekidis:2003:EFM}.
The main idea is to use short
bursts of appropriately initialized simulations to estimate the time
derivative of the quantities of interest and then use polynomial extrapolation
to jump forward in time
\cite{Gear:2002:CIB,Erban:2006:CAE}.
When projective integration is applied to
deterministic problems (governed by systems of differential equations),
one can show that it might significantly accelerate the computation of
time evolution for systems with large gaps in
their eigenvalue spectrum \cite{Gear:2003:PMS}.
By wrapping the same algorithm around an inner atomistic and/or stochastic
simulator, one can similarly accelerate coarse-grained computations
\cite{Gear:2001:PIM,Erban:2006:CAE,Gear:2002:CIB}.

Many problems possess additional continuous symmetries
\cite{Golubitsky:2003:TSP,LeMesurier:1988:FMF,Chen:1995:NNG,Beyn:2004:FSE}
which can give rise to solutions which
are traveling, exploding, collapsing or rotating in the domain of interest.
In principle, projective integration might be applied to such systems
as well.
However, we can improve the efficiency of the method
by taking the underlying symmetry into account.
The key idea is to perform the projective integration in a
``co-evolving" frame \cite{Rowley:2000:REK}.

Projective integration in a co-traveling frame is applied to the
Nagumo equation \cite{Miura:1982:ACF,Murray:2002:MB}, a well-studied
system with translational invariance and traveling solutions.
Projecting in a dynamically renormalized frame is similarly
applied to systems characterized by scale invariance.
The scale invariant
system we study in this paper is one-dimensional diffusion.
In both applications the existence of continuous symmetries
is exploited.
We apply modified projective integration protocols that are
implemented in a dynamically co-evolving
frame.
We demonstrate that this modification improves computational
accuracy, allowing for large projective steps.

We also illustrate the coarse-grained version of projective integration in
a co-traveling frame, for a Stochastic Simulation Algorithm
(SSA) \cite{Gillespie:1977:ESS} implementation of the Nagumo kinetics.
The density of reactant
particles progressively forms a traveling wave front moving with a constant shape
and velocity.
In a second application we study one-dimensional diffusion
simulated by a large ensemble of random walkers.
The macroscopic behavior, described by the cumulative density function (CDF)
of the particle positions, features scale invariance.
Projective integration is appropriately modified to exploit
this scale invariant character of the macroscopic behavior and
increase the accuracy of computations, again allowing for relatively large
projective  time steps.

The paper is organized as follows.
In Section \ref{secprojcoarseproj}
we briefly summarize projective integration techniques and
discuss the general ideas of equation-free techniques
\cite{Kevrekidis:2003:EFM}, a computational framework
wrapped around microscopic (e.g. kinetic Monte Carlo)
simulators.
We then present our projective and coarse projective integration scheme in a
co-evolving frame and its application to
problems with translational invariance (Section \ref{transl_invariance})
as well as scale invariance (Section \ref{scale_invariance}).
In Section \ref{FHNsim} we illustrate the efficiency of the method
in a co-traveling frame for the Nagumo equation.
We also describe coarse projective integration
for a kinetic Monte Carlo simulation of a reaction-diffusion system based on
Nagumo kinetics.
In Section \ref{diffusionSIM} we present results from the
application of projective integration to the scale invariant diffusion
system, both for the macroscopic diffusion equation and for a
random-walker model in one spatial dimension.
In Section \ref{generalapproach} we propose a more general approach
that can handle systems evolving in space and scale with an {\em asymptotically}
invariant form, and summarize our work in Section \ref{secdiscussion}.

\section{Projective and Coarse Projective Integration}

\label{secprojcoarseproj}

\numberwithin{equation}{section}
Consider a system described by either a suitable macroscopic
evolution equation or a stochastic, individual-based model, and let
$M(t)$ be the macroscopic observable, for which a closed macroscopic
evolution equation exists.
Depending on the problem, $M(t)$ can be a single scalar, a vector or
a point in a suitable infinite-dimensional Banach space (e.g. a
function of physical space).
In our illustrative numerical examples the macroscopic observable
will be a (discretized) field of the density of individuals.

The parameters of the {\it forward Euler projective integration scheme}
are two time constants, $\Delta t$ and $T$.
Given the value of the macroscopic observable at time $t$, a
suitable ``inner'' timestepper (e.g. the stochastic simulator) is
used to compute the system evolution until time $t+\Delta t$. Using
the values of $M$ in the interval $(t, t+\Delta t)$ we estimate the
time derivative of $M$ and use it to estimate (project) the value of
the macroscopic observable $M(t+\Delta t+T)$ using a Taylor
expansion:

\begin{equation}
\label{Taylor_expansion}
M(t+\Delta t+T) \approx M(t+ \Delta t)
+ T \frac{\partial{M}}{\partial{t}} \Big{|}_{ (t+ \Delta t)}.
\end{equation}

\noindent
Hence, we compute $M(t+\Delta t+T)$ from the value of $M(t)$ by
running the inner integrator for time $\Delta t$ only.
Other, more sophisticated projective integration schemes can be
readily constructed \cite{Gear:2002:CIB,Gear:2003:TPS,Lee:2005:SOPI}.

\begin{figure}
\begin{center}
\includegraphics[width=0.7\linewidth]{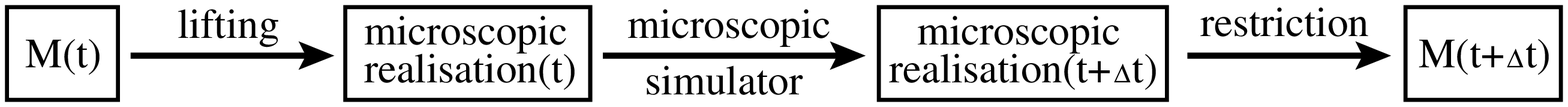}
\end{center}
\caption{{\it Schematic of a coarse time stepper.}}
\label{microtimestepper}
\end{figure}

If the evolution equation for $M(t)$ is explicitly available, it is
straightforward to compute $M(t+\Delta t)$ from $M(t)$ using a
suitable discretization of this available evolution equation.
However, if the only information for the time evolution of the
system comes from an individual-based, stochastic model, then we
have to use the idea of the coarse timestepper
\cite{Kevrekidis:2003:EFM} as illustrated schematically in Fig.\ref{microtimestepper}.
Given a macroscopic variable $M(t)$, we construct consistent
microscopic initial conditions; we call this the {\em lifting}
procedure.
Next, we evolve the system using the microscopic simulator (e.g.
kinetic Monte Carlo) for time $\Delta t$.
Now we compute $M(t+ \alpha_i\Delta t)$ from the microscopic data
for various instances $0 < \alpha_1 < \cdots < \alpha_K=1$ (the {\em
restriction} step).
Having computed $M$ at these instances, we use them to estimate the
time derivative of $M$ and use (\ref{Taylor_expansion}) as in the
deterministic case.
The resulting method is {\em coarse} projective integration.
\\

\noindent

Some computational gain from projective or coarse projective
integration can be expected provided we can choose $\Delta t \ll T$.
A relatively large extrapolation step may save substantial
computational time, considering the computational demands of a
particle-level simulator.
On the other hand, large steps can lead to low accuracy in
simulating the dynamics of the system, or even cause numerical
instabilities.
We will discuss how one can obtain increased accuracy in
projectively integrating systems with solutions evolving
along continuous symmetry groups.
The key idea is to evolve the solution (macroscopic observable) in a
coordinate frame which tracks the evolution across space (for problems with
traveling solutions) and across scales (in problems with self-similar solutions).

\subsection{Projective integration in a co-evolving frame}

\label{modProjective}

In many cases of interest, the long time macroscopic dynamics
do not involve stationary solutions but rather traveling,
rotating, or
scale invariant (e.g. self-similar) solutions \cite{Barenblatt:1996:SSI}.
Accuracy concerns in the direct application of projective
integration to problems with such solutions
\cite{Rowley:2003:RRS,Beyn:2004:FSE} limit the projective time step
$T$.
Consider a traveling wave solution for a problem with translational
invariance: it is natural to study its evolution in a co-traveling
frame, where the solution asymptotically appears stationary (the
traveling has been factored out).
In the same sense, it is natural to study self-similar solutions in
a dynamically renormalized frame, where the {\em scale} evolution of
the solution has been factored out.
Recently, a template based approach has been developed for the
investigation of problems with translational invariance
\cite{Rowley:2000:REK} (see also \cite{Beyn:2004:FSE}) and has been extended to the study of
self-similar solutions
\cite{Rowley:2003:RRS,Aronson:2001:GFL,Siettos:2003:FRR}.
If the description of the macroscopic dynamics involves scale
invariant partial differential equations (PDEs), template conditions
can be applied to derive equations describing the evolution in a
dynamically renormalized framework.
The steady state of the renormalized equations correspond to
self-similar solutions of the original problem and the similarity
exponents can also be conveniently computed \cite{Aronson:2001:GFL}.
This dynamic renormalization concept can also be applied to
multiscale system models where an explicit formulation for the
macroscopic evolution equation is not available
\cite{Chen:2004:MDC,Chen:2006:EFD,Szell:2005:CCC,Kessler:2006:EDR,Erban:2006:DPI}.

In this paper, our goal is to study coarse projective integration
for such multiscale atomistic and/or stochastic problem models.
To explain the idea of the co-evolving frame for such problems, it
is easier to start with a deterministic example.
We consider the PDE written in the following form
\begin{equation}
\label{generalPDE} \frac{\partial M}{\partial{t}}= \operL_{x}
\left(M \right).
\end{equation}
Here, $M \equiv M(x,t) \in \seB$ and $\operL_{x}: \seB \to \seB$
where $\seB$ is a suitable Banach space of functions mapping $\er$ to $\er$
and the subscript $x$ denotes the independent space variable.
We define the
{\it shift operator} $\caS_C: \seB \to \seB $ and the {\it rescaling
operator} $\caR_{A,B} : \seB \to \seB$ by
\begin{equation}
\caS_C (f) : x \to f(x+C)
\qquad
\mbox{and}
\qquad
\caR_{A,B} (f) : x \to B f \left( \frac{x}{A} \right)
\label{operRS}
\end{equation}
for any $A, B >0$ and $C \in \er$.
We distinguish two cases -- projective integration in a co-traveling
frame in Section \ref{transl_invariance} and projective integration
in a frame which scales with the solution in Section
\ref{scale_invariance}.

\begin{figure}
\begin{center}
\includegraphics[width=1\linewidth]{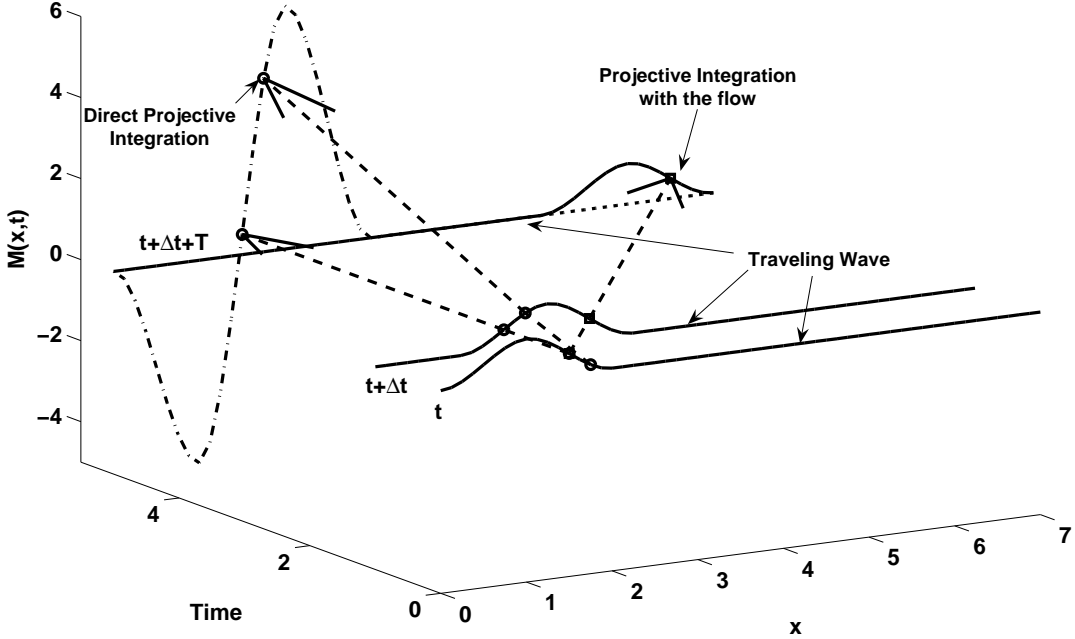}
\end{center}
\caption{(Schematic) Direct projective integration fails to produce the correct traveling shape at 
time $t+\Delta t +T$. Projective integration in a co-evolving frame gives results with higher accuracy
for the same time step.}
\label{wrong_projection}
\end{figure}

The appropriateness of a co-evolving frame for projective
integration is schematically illustrated in Fig.\ref{wrong_projection} for
a constant shape traveling wave.
Direct projective integration uses the computed wave at different time instances
to estimate its time derivative. 
For the instances shown in the Figure, projection according to 
(\ref{Taylor_expansion}) produces manifestly wrong results for large $T$;
On the other hand, projection with the same data in a
co-evolving frame gives results with much higher accuracy for the same time step $T$
(i.e. for the same computational cost).
This is because the time derivative is much smaller (here practically zero)
in the co-evolving frame.
\subsection{Systems with translational invariance}
\label{transl_invariance}

Let the differential operator $\operL_{x}$ in (\ref{generalPDE})
satisfy the translational invariance property, i.e. the following
relation holds for every $C \in \er$:
\begin{equation}
\operL_{x} \caS_C = \caS_C \operL_{x}.
\label{trav_invariance}
\end{equation}
Let $M(x,t) \in \seB$ be the solution of (\ref{generalPDE}) and let
$C(t)$ be a differentiable function of time.
We define
\begin{equation}
\label{shift2} \widehat{M}(x,t) = \caS_{C(t)} M(x,t), \qquad
\mbox{which means that} \qquad M(x,t) = \caS_{-C(t)}
\widehat{M}(x,t).
\end{equation}
Using (\ref{shift2}) and (\ref{generalPDE}), we obtain
\begin{equation}
\label{cotraveling}
\frac{\partial{\widehat{M}}}{\partial{t}}=
\operL_{x} \widehat{M}
+
\frac{\D C}{\D t} \frac{\partial{\widehat{M}}}{\partial{x}}.
\end{equation}
If $C(t)$ is given, then solving (\ref{cotraveling}) provides the
same information as (\ref{generalPDE}).
We have the freedom to choose $C(t)$; we will do it so as to
naturally take into account the ``traveling component'' of the
solution.
To find an appropriate shift $C(t)$, we impose an additional
algebraic constraint (template condition)
\cite{Rowley:2000:REK,Beyn:2004:FSE}.
The purpose of this template condition is to determine a shift,
$C(t)$, such that $S_{C(t)}M$ is as ``independent of $t$ as
possible.''
If $M(x,t)$ were a constant shape traveling wave, then $C(t)$ would
be the change in wave position with time and $S_{C(t)}M$ would
be stationary.
Hence we need a way to measure how far the wave has moved so we can
shift it back by that amount.
One can construct suitable templates in many ways.
A seemingly natural way is to ask for the shift that minimizes some
norm of the difference between the shifted wave and some fixed
waveform (the ``template'', $\widetilde{T}(x)$ ):
\begin{equation}
\frac{\partial}{\partial{C}}  ||\caS_{C}M(x)-\widetilde{T}(x)|| =0.
\end{equation}
That fixed waveform could be $M(0,x)$, something believed to
approximate the final solution, or, in principle, anything else.
An alternative is to use the centroid of the absolute value of the
wave or some other characteristic that identifies ``where the wave
is'' and apply a shift to bring this feature to a constant position in
space.
For a single-humped wave, one might consider using the location of
the wave maximum. However, if during  a transient the wave develops a second
maximum this would clearly fail.
The centroid of the absolute value is unique and easy to compute.
If the wave is positive (or of constant sign), such as a density
measure, the centroid has the advantage of being a linear 
of the wave shape that will be preserved under projective
integration.
We will formally write the template condition as the
algebraic equation:
\begin{equation}
\label{template_trans} \hbar(\widehat{M},\widetilde{T}) = 
\hbar(S_{C(t)}M,\widetilde{T}) = 0,
\end{equation}
where $\hbar$ is a functional mapping $\mathbb{B} \times \mathbb{B}$ to $\er$.
This, together with (\ref{cotraveling}) describe the dynamics of
the shifted solution $\widehat{M}(x,t)$ as well as the dynamics of its shift,
$\D C(t)/\D t$.
Such template conditions arise naturally in the computation of limit
cycle solutions in autonomous dynamical systems, where they are
often also called ``pinning" conditions \cite{Doedel:1986:AUTO,Beyn:2004:FSE}.

If the operator $\operL_{x}$ is available explicitly, we can use
(\ref{cotraveling}) to estimate the time derivative of $\widehat{M}$
at a given time $t$ and (\ref{Taylor_expansion}) can be
used to make the extrapolation in time for simple projective forward
Euler.
To complete the projective algorithm in the co-evolving frame, we
also must specify how the shift, $C(t)$, evolves during a projective
time step.
As in (\ref{Taylor_expansion}) we approximate the shift evolution by
\begin{equation}
\label{shift_evol}
C(t + \Delta t + T)
\approx
C(t+\Delta t) + T \frac{\D C}{\D t}\Bigg{|}_{(t+\Delta t)}.
\end{equation}
Finally, the {\em unshifted} projected solution $M(t+\Delta t + T)$
is computed, if required, from
\begin{equation}
\label{shift}
M = \caS_{-C(t+ \Delta t + T)} \widehat{M}.
\end{equation}

\subsection{Systems with scale invariance}
\label{scale_invariance}

Consider a scale invariant problem where the differential operator
$\operL$ satisfies the property
\begin{equation}
\operL \caR_{A,B} = A^aB^{b-1} \caR_{A,B} \operL,
\end{equation}
 i.e. there exist constants $a$ and $b$ such
that the above relation holds for every $A,B>0$.
Equivalently
\begin{equation}
\label{sc_invariance}
\operL_{x}\left( BM \left( \frac{x}{A} \right) \right)=B^{b}A^{a}\operL_{y} \left( M(y) \right)
\quad \mbox{where} \quad y=\frac{x}{A}
\end{equation}
where $\operL_{x}$ and $\operL_{y}$ denotes the action of the operator
$\operL$ on the coordinates $x$ and $y$ respectively.
Note that the system must {\em also} satisfy the translational
invariance property (\ref{trav_invariance}); however, for
simplicity, we will assume that we do not have traveling
solutions, but concentrate on self-similar ones. The combination
is considered in Section \ref{generalapproach}.

We study solutions $M(x,t)$ of (\ref{generalPDE}); choosing
scaling factors $A$ (for space) and $B$ (for the solution amplitude)
as well as a reparametrization of time $\tau(t)$ leads to the study
of the {\it rescaled} solutions $\widehat{M}(y,\tau)$
(see \cite{Aronson:2001:GFL})
\begin{equation}
\label{gen_scale} M(x,t)=B(\tau) \widehat{M} \left(
\frac{x}{A(\tau)} , \tau (t) \right).
\end{equation}
Equation (\ref{generalPDE}) can be re-written in such a dynamically
renormalized form as follows \cite{Rowley:2003:RRS,Aronson:2001:GFL}
\begin{equation}
\label{coevolve}  \left(
\frac{1}{B}\frac{\D {B}}{\D {\tau}} \widehat{M}-
\frac{1}{A} \frac{\D {A}}{\D {\tau}}  y
\frac{\partial{\widehat{M}}}{\partial{y}}+\frac{\partial{\widehat{M}}}{\partial{\tau}}
\right)\frac{\D {\tau}}{\D {t}}= A^{a}B^{b-1}\operL_{y}
\left( \widehat{M}(y) \right).
\end{equation}
Motivated by the search for self-similar solutions we select the
time reparametrization $\tau (t)$ as
\begin{equation}
\label{tauoft}
\frac{\D {\tau}}{\D {t}}=A^{a}B^{b-1}
\end{equation}
which leads to
\begin{equation}
\label{coevolve_fin}
\frac{\partial{\widehat{M}}}{\partial{\tau}}=
\operL_{y} \left( \widehat{M}(y) \right) -\frac{1}{B}\frac{\D {B}}{\D {\tau}} \widehat{M}
+\frac{1}{A}\frac{\D {A}}{\D {\tau}}y
\frac{\partial{\widehat{M}}}{\partial{y}}.
\end{equation}
For self-similar solutions
$\frac{1}{A}\frac{\D {A}}{\D {\tau}}$ as well as
$\frac{1}{B}\frac{\D {B}}{\D {\tau}}$ are constants whose
particular values depend on $\widehat{M}$.
In our renormalization algorithm $\widehat{M}$ is selected
by our choice of template condition(s). 

\begin{figure}[h]
\begin{center}
\includegraphics[width=0.6\linewidth]{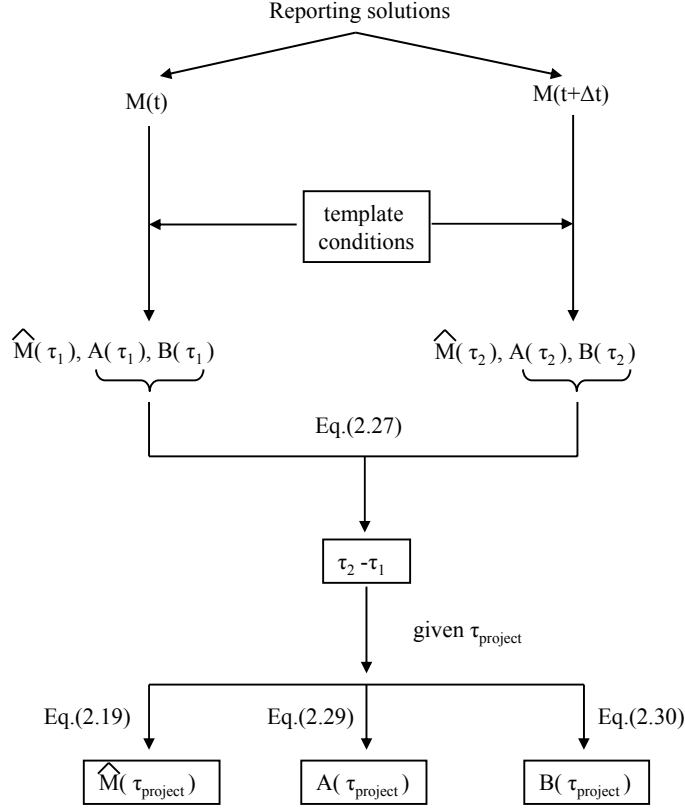}
\end{center}
\caption{{\it Schematic description of template based time projection for scale invariant
problems.}}
\label{flowtau}
\end{figure}

The main idea of projective integration in a co-evolving frame is
to factor out the scale evolution, so as to obtain a (rescaled)
solution that evolves more slowly.
As in the traveling case, we exploit the template-based approach in
order to compute solutions which are ``as scale invariant as
possible''.
A schematic description of the projective integration algorithm
to scale invariant problems, is shown in Fig.\ref{flowtau}.

To determine the evolution of the scale parameters $A$ and $B$ we
need to apply template conditions that control the spatial extent of
the solution ($A$) as well as its amplitude ($B$).
We could, for example, minimize the distance between the {\em
rescaled} current solution and a template function to determine both
$A$ and $B$.
Alternatively we could determine the amplitude by maintaining the
constancy of some norm of the solution - the $L^{1}$-norm for a
positive function would simply maintain the total mass - while for
the spatial extent we could keep a moment of a positive measure of
the solution, such as
$\int_{-\infty}^\infty|\widehat{M}(x,t)|x^2 \D x$, constant (assuming
that this integral is well defined).
In general, we need two independent conditions to determine the two
scale parameters; these will take the form of two algebraic
equations which are used along with (\ref{coevolve_fin})
\cite{Aronson:2001:GFL,Rowley:2003:RRS,Chen:2004:MDC} to evolve the
dynamically renormalized problem.
These two algebraic equations take the form
\begin{equation}
\hbar_{A}(\widehat{M}(y),\widetilde{T}_{1})=0 \quad \Longleftrightarrow \quad
\hbar_{A}\left( \frac{1}{B}M(Ay),\widetilde{T}_{1}\right)=0
\end{equation}
\begin{equation}
\hbar_{B}(\widehat{M}(y),\widetilde{T}_{2})=0 \quad \Longleftrightarrow \quad
\hbar_{B}\left( \frac{1}{B}M(Ay),\widetilde{T}_{2} \right)=0
\end{equation}
where $\widetilde{T}_{1}$, $\widetilde{T}_{2}$ are template functions.
While there is considerable freedom in the choice of these algebraic
conditions, it is important that they yield a {\em unique} and
computationally simple solution for the scaling factors $A$ and $B$.

We can apply projective integration to the rescaled solution $\widehat{M}$ 
in the original time frame $t$, or in the rescaled time frame $\tau$ chosen above
(or, for that matter, in any other convenient time variable).
Using $\tau$, projective forward Euler is:
\begin{equation}
\label{Taylor_tau}
\widehat{M}(\tau_{project})=
\widehat{M}(\tau_{2})+(\tau_{project}-\tau_{2})
\frac{\partial \widehat{M}}{\partial{\tau}}
\Bigg{|}_{\tau_{2}} \approx
\widehat{M}(\tau_{2})+(\tau_{project}-\tau_{2})
\frac{\widehat{M}(\tau_{2})-\widehat{M}(\tau_{1}) }{\tau_{2}-\tau_{1}}
\end{equation}
where the rescaled times $\tau_{1},\tau_{2}$ and $\tau_{project}$
correspond to times $t$, $t+ \Delta t$ and $t+\Delta t +T$ respectively.
In order to approximate numerically the right hand side of (\ref{Taylor_tau})
we must determine the relation between time $t$ and the rescaled time $\tau$.

\begin{figure}
 \begin{center}
 \begin{tabular} {cc}
 \includegraphics[width=0.45\linewidth]{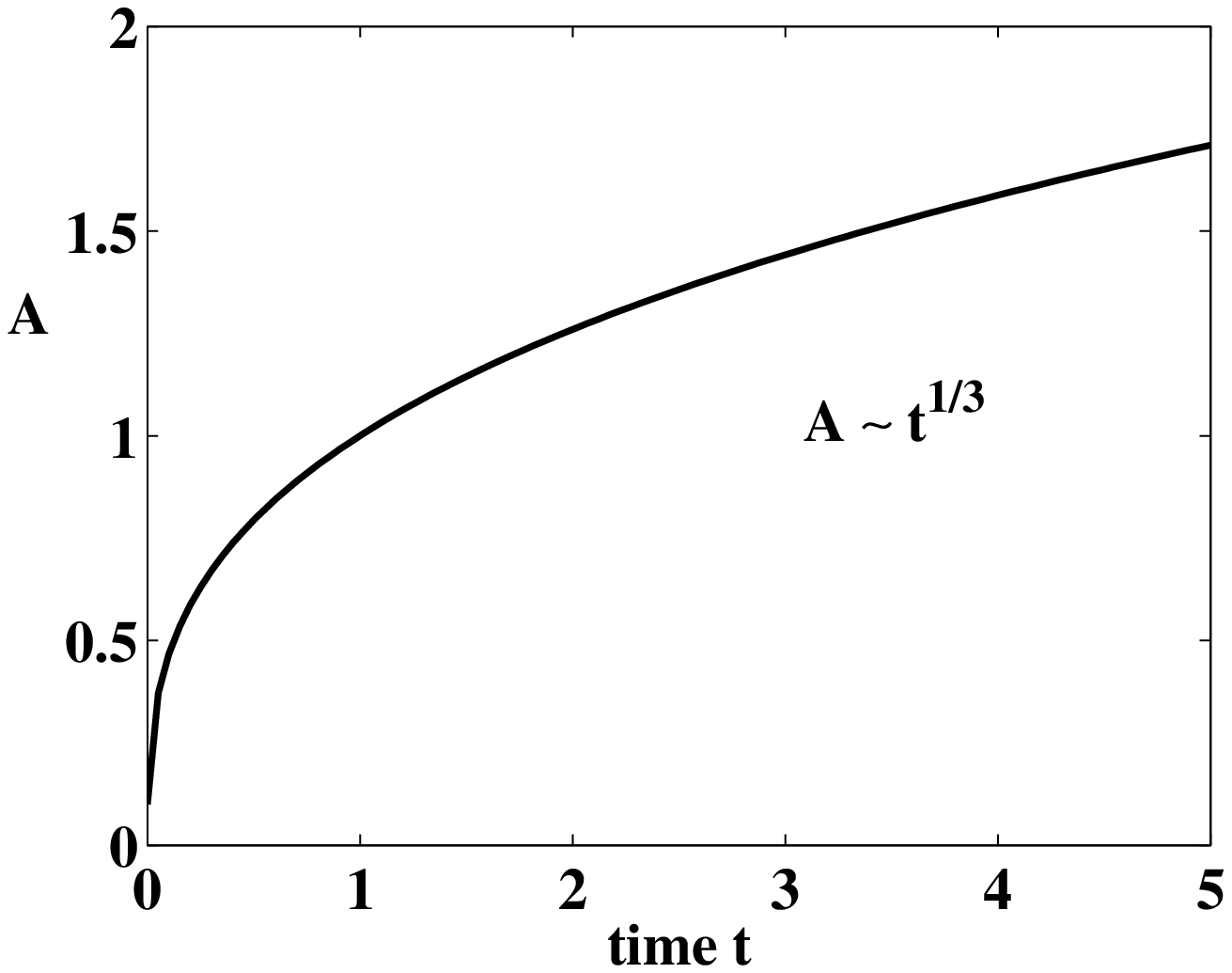} &
\includegraphics[width=0.45\linewidth]{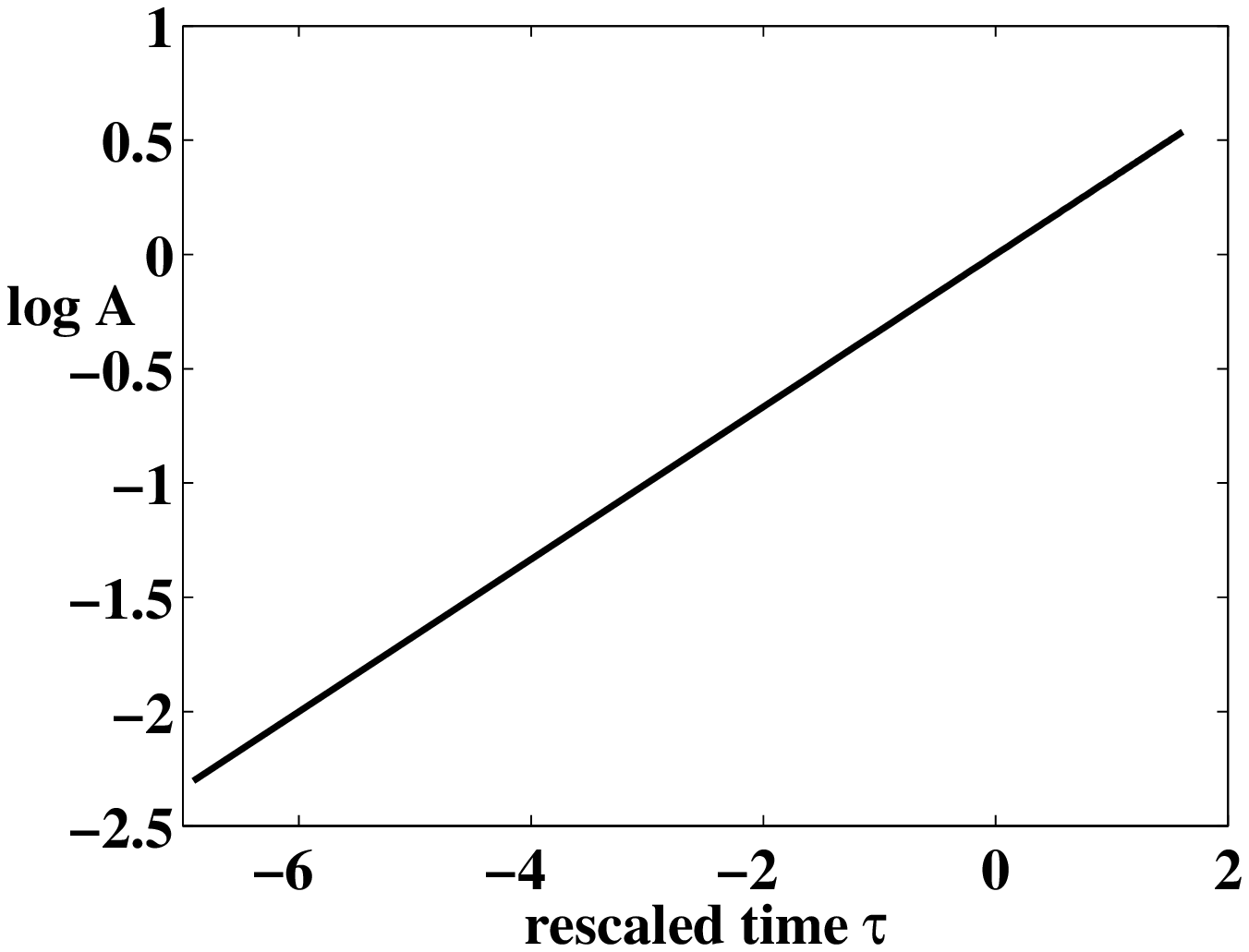} \\
    (a) & (b)
    \end{tabular}
    
  \end{center}
 \caption{ (a) Evolution of scale factor $A$ as a function of time $t$  
for the self-similar solution of 1D porous medium equation. 
The similarity exponent $\gamma$ (see text for details) is equal to
$\gamma=1/3$. The blow-up time is $t^{*}=0$. (b) Linear evolution 
of $\log{A}$ with respect to rescaled time $\tau$ ($\tau \sim \log{t}$). }
 \label{ksiexplain}
 \end{figure}

We will assume that during the projection step the parameters
\begin{equation}
\xi_{A}=\frac{1}{A} \frac{\D{A}}{\D{\tau}},\quad
\xi_{B}=\frac{1}{B} \frac{\D{B}}{\D{\tau}}
\end{equation}
remain constant (this is true for self-similar solutions, and
is analogous to assuming that the velocity
is constant during a projection step for traveling problems).
The evolution of the scale factors $A,B$ for a self-similar
problem is of the form \cite{Aronson:2001:GFL}:
\begin{equation}
A(t)\sim |t-t^{*}|^{\gamma}, \quad B(t)\sim |t-t^{*}|^{\delta}
\end{equation}
where $t^{*}$ is an appropriate (positive or negative) blow-up time, and
$\gamma, \delta$ are the similarity exponents. 
A typical example is the 
1D Barenblatt solution, the self-similar solution to the porous medium 
equation $u_{t}=(u^{2})_{xx}$ \cite{Aronson:1986:PME,Aronson:2001:GFL}. In this case, the similarity 
exponents are $\gamma=1/3$ and $\delta=-1/3$. If we consider the case, where the 
blow-up time $t^{*}=0$ (the initial datum is a Dirac mass
at the origin), then the scale factor $A$ can evolve as depicted in 
Fig.\ref{ksiexplain}(a). In Fig.\ref{ksiexplain}(b), one can see the
{\it linear} evolution of $\log{A}$ with respect to the rescaled time
$\tau$. It can be shown \cite{Aronson:2001:GFL} that the relation between
rescaled time $\tau$ and time $t$ has the form $\tau \sim \log{t}$;
$\log(B)$ behaves similarly.
We thus expect better accuracy when the projective scheme 
is based on exponential growth of the scale factors $A$ and $B$.

During the projective step the evolution of $A$ and $B$ is described by:
\begin{equation}
\label{Aevol}
A(\tau)=A(\tau_{1})\exp(\xi_{A} (\tau-\tau_{1}))
\end{equation}
\begin{equation}
\label{Bevol}
B(\tau)=B(\tau_{1})\exp(\xi_{B} (\tau-\tau_{1})).
\end{equation}
The values $A(t) \equiv A(\tau_1), B(t) \equiv B(\tau_1),A(t+
\Delta t) \equiv A(\tau_2), B(t+\Delta t) \equiv B(\tau_2)$ obtained through the application of the 
template conditions are used in the relation
\begin{equation}
\label{Astar}
\xi_{A}(\tau_{2}-\tau_{1})=\log\frac{A(\tau_{2})}{A(\tau_{1})}=
\log\frac{A(t+\Delta t)}{A(t)}=A^{*}
\end{equation}
\begin{equation}
\label{Bstar}
\xi_{B}(\tau_{2}-\tau_{1})=\log\frac{B(\tau_{2})}{B(\tau_{1})}=
\log\frac{B(t+\Delta t)}{B(t)}=B^{*}.
\end{equation}
We can now derive an expression between the rescaled time step $\tau_{2}-\tau_{1}$
and the time step $\Delta t$. Namely from (\ref{tauoft}) and (\ref{Aevol}) -- (\ref{Bevol}), 
we get
\begin{equation}
\label{timestep}
\frac{A(\tau_{1})^{-a} B(\tau_{1})^{1-b}}{-a\xi_{A}+(1-b)\xi_{B}}
\Big\{
\exp \big[ -a\xi_{A}(\tau_{2}-\tau_{1})+(1-b)\xi_{B}(\tau_{2}-\tau_{1}) \big] -1
\Big\}
=\Delta t.
\end{equation}
It is straightforward to evaluate  the parameters $\xi_{A}$ and $\xi_{B}$
from (\ref{Astar}) -- (\ref{Bstar}) to obtain 
\begin{equation}
\label{timestep}
\frac{A(\tau_{1})^{-a} B(\tau_{1})^{1-b}}{-aA^{*}+(1-b)B^{*}}
\Big\{
\exp \big[ -aA^{*}+(1-b)B^{*} \big] -1
\Big\}
(\tau_{2}-\tau_{1})=\Delta t.
\end{equation}
Then we compute the rescaled projection time $\tau_{project}$ from
\begin{equation}
\label{tauproject}
\frac{A(\tau_{1})^{-a}B(\tau_{1})^{1-b}}{-a\xi_{A}+(1-b)\xi_{B}}
\Big\{ \exp \big[ (-a\xi_{A}+(1-b)\xi_{B}) (\tau_{project}-\tau_{1}) \big]
-1 \Big\} = T+ \Delta t.
\end{equation}
Finally, we can also obtain the projections of 
the scale factors $A,B$ from (\ref{Aevol}),(\ref{Bevol}):
\begin{equation}
\label{Aproject}
A_{project}=A(\tau_{project})=A(\tau_{1})\exp\big[\xi_{A} (\tau_{project}-\tau_{1})\big]
\end{equation}
\begin{equation}
\label{Bproject}
B_{project}=B(\tau_{project})=B(\tau_{1})\exp\big[\xi_{B} (\tau_{project}-\tau_{1})\big]
\end{equation}
and, if desirable, recover the projection of full solution $M(t+\Delta t +T)$ from the
rescaling relation:
\begin{equation}
M(t+ \Delta t +T)=B_{project}\widehat{M}\left( \frac{x}{A_{project}},\tau_{project} \right).
\end{equation}

\section{Systems with translational invariance - A reaction-diffusion problem}

\label{FHNsim}

In this section we demonstrate the efficiency of the proposed projective integration
scheme in a co-evolving frame for a reaction-diffusion system with translational
invariance. 
Our stochastic model is motivated by the Nagumo equation
\cite{Miura:1982:ACF,Murray:2002:MB,Beyn:2004:FSE},
\begin{equation}
\label{FHN}
\frac{\partial{u}}{\partial{t}}=D \ \frac{\partial ^{2}u}{\partial{x}^{2}}+u(1-u)(u-\alpha)
\end{equation}
where $u$ denotes the reactant's concentration, $\alpha$ is a
kinetic parameter and $D$ is the 
diffusion coefficient. We use $D=1$ in what follows. 
We consider a large (effectively infinite)
domain with zero flux (Neummann) boundary conditions.
The Nagumo equation is a
a well known example of a parabolic system that exhibits traveling waves and
has an explicit wave solution $u(x,t)=\hat{u}(x-c)$ given by:
\begin{equation}
\label{theorNag}
\hat{u}(x)= \left[{1+\exp\left(-\frac{x}{\sqrt{2}}\right)} \right] ^{-1} \quad \mbox{,} \quad
 \frac{\D c}{\D t}=-\sqrt{2}\left(\frac{1}{2}-\alpha\right).
\end{equation}
Motivated by the Nagumo kinetics we construct a particle-based simulator, where
a set of chemical reaction steps as well as diffusion steps are incorporated.
We consider the following set of reactions:

\begin{equation}
{\mbox{ \raise 0.851 mm \hbox{$2N + H$}}}
\quad \;
\mathop{\stackrel{\displaystyle\longrightarrow}\longleftarrow}^{k_{1}}_{k_{-1}}
\quad
{\mbox{ \raise 0.851 mm \hbox{$3N$}}}
\label{firstreact}
\end{equation}

\begin{equation}
\label{secondreact}
N \stackrel{k_{2}} \longrightarrow \varnothing .
\end{equation}

The reaction rate constant for the production of reactant $N$ is $k_{1}=1+\alpha$, while the consumption
rate constants are respectively $k_{-1}=1$ and $k_{2}=\alpha$.
The concentration of reactant $H$ is assumed to remain essentially constant 
and equal to $1$ ($H=1$).

A standard way to simulate a spatially homogeneous chemical system
is the Gillespie SSA \cite{Gillespie:1977:ESS}.
At each time step of the algorithm a pair of random numbers is generated in order 
to answer two essential questions:
when will the next event -- chemical reaction -- occur and which reaction will
it be?
We incorporate in our system the effect of spatial diffusion too;
$N$ diffuses with a diffusion coefficient, $D$.

The generalisation of Gillespie ideas to spatially
distributed systems can be found in e.g.
\cite{Stundzia:1996:SSC,Isaacson:2006:IDC}.
Here, diffusion is treated as another set of ``reaction steps'' in the system.
The domain of interest is discretized into $J$ lattice sites with constant
distance $h$ between them.
We denote by $N_i$  the number of respective molecules
at lattice site $i$.
This means that we describe the state of the stochastic
reaction-diffusion system by a $J$-dimensional vector
$\boldN = \big[ N_1, N_2 \dots, N_{J} \big],$
and the following reactions at each time step are considered
\begin{equation}
\left.\begin{array} {l}
{\mbox{ \raise 2.551 mm \hbox{$2N_{i} + H$}}}
\quad \;
\stackrel{\stackrel{{ k_{1}}}{\displaystyle \longrightarrow}}{\stackrel{\displaystyle \longleftarrow}{_{k_{-1}}}}
\quad
{\mbox{ \raise 2.551 mm \hbox{$3N_{i}$}}}\\
\quad N_{i} \quad \stackrel{k_{2}} \longrightarrow \qquad \varnothing
\end{array}
\right\}~~i=1, \dots, J, \label{reactboxi}
\end{equation}
\vspace{-0.5cm}

\begin{equation}
N_i~\stackrel{d}{\longrightarrow} N_{i+1}, ~~~~~~~~i=1, \dots, J-1,
\label{diffip1}
\end{equation}
\begin{equation}
N_i~\stackrel{d}{\longrightarrow} N_{i-1},~~~~~~~~i=2, \dots, J.
\label{diffin1}
\end{equation}

The set of reactions (\ref{reactboxi}) implies that the reaction
mechanism (\ref{firstreact}) -- (\ref{secondreact})
is implemented at each lattice site of the domain. 
Moreover,
diffusion is introduced as a set of new reactions (\ref{diffip1}) -- (\ref{diffin1}), whose
transition rates are denoted by
$d$. 
The transition rates for $d$ are connected to the macroscopic diffusion
coefficient, $D$, which at a certain limit ($h\ll1$), is given by
the formula $d = D/h^2 $. 
The augmented set of reactions
(\ref{reactboxi}) -- (\ref{diffin1}), together with
suitable boundary conditions, can thus be simulated using  Gillespie SSA.

We will start by
projectively integrating the Nagumo partial differential equation,
and then we will illustrate
the {\it coarse} variant of the method for the particle-based implementation of
the scheme (Gillespie algorithm).

\subsection{Nagumo Equation  - PDE description}
We first consider the deterministic description (PDE) of the Nagumo problem
and illustrate the accuracy improvement to the projective method from
operating in a co-evolving frame.
In our numerical computations $\alpha=0.01$, so that velocity of the
traveling wave is $\D c / \D t \approx -0.693$ according to (\ref{theorNag}).
The (long) one-dimensional domain $[-30,30]$ is discretized  into $601$
equidistant nodes (i.e., the distance between two successive nodes
is $\delta x=0.1$). 
The spatial partial derivatives are approximated with central finite differences
and the applied boundary conditions are of Neumnann type. 
The initial condition is:
\begin{equation}
u(x,0)=\left \{ \begin{array} {ccc}
0& \textrm{ for } & -30<x\le 0, \\
x/10& \textrm{ for } & 0<x\le 10, \\
1 & \textrm{ for } & 10<x \le 30.
\end{array}
\right.
\end{equation}
The time step of the inner integrator (here a simple forward Euler) is $\delta t=10^{-4}$
to satisfy the stability criterion $2 \delta t< \delta x^{2}$.

\begin{figure}
\begin{center}
\begin{tabular}{cc}
\includegraphics[width=0.45\linewidth]{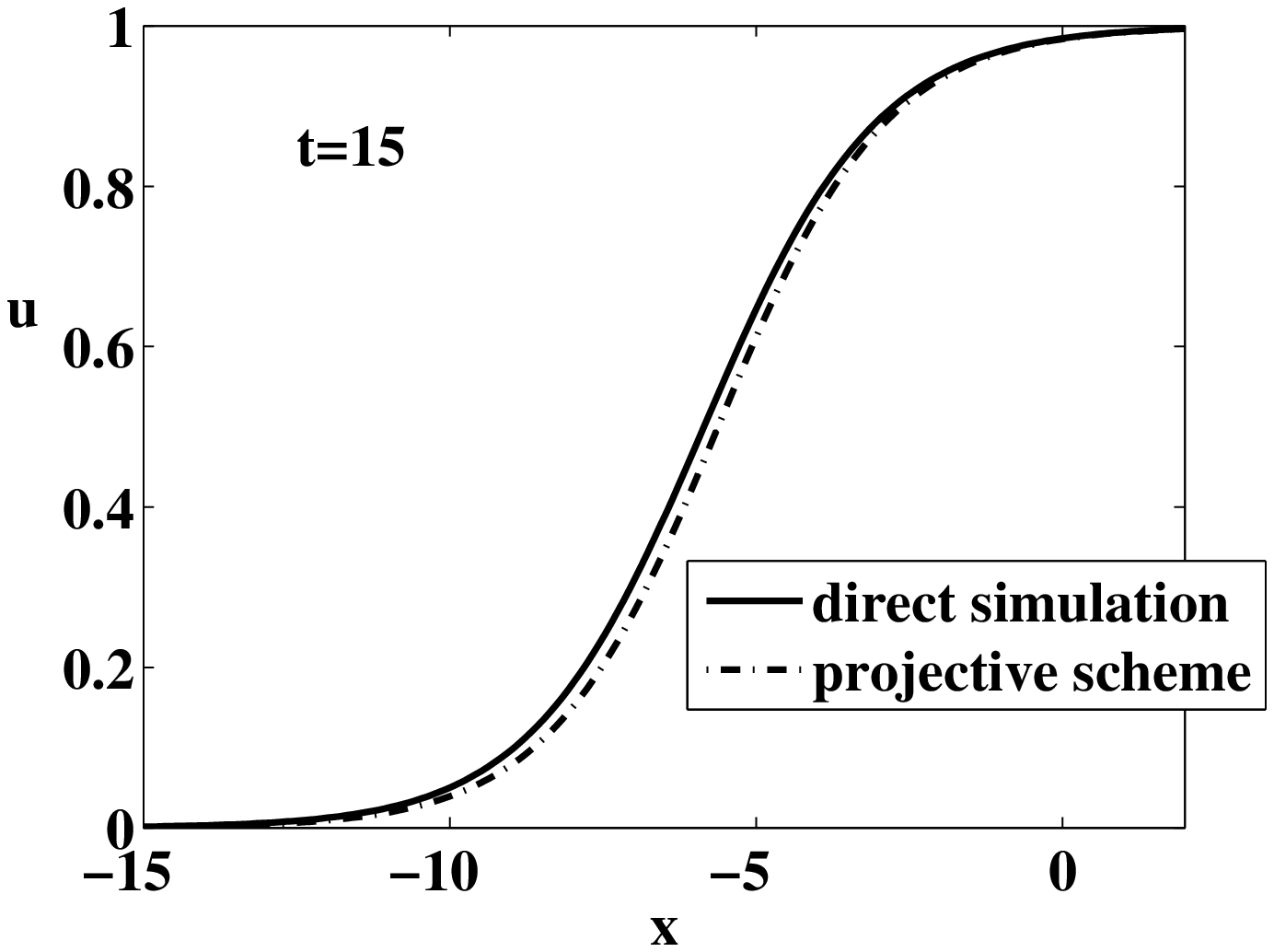} &
\includegraphics[width=0.45\linewidth]{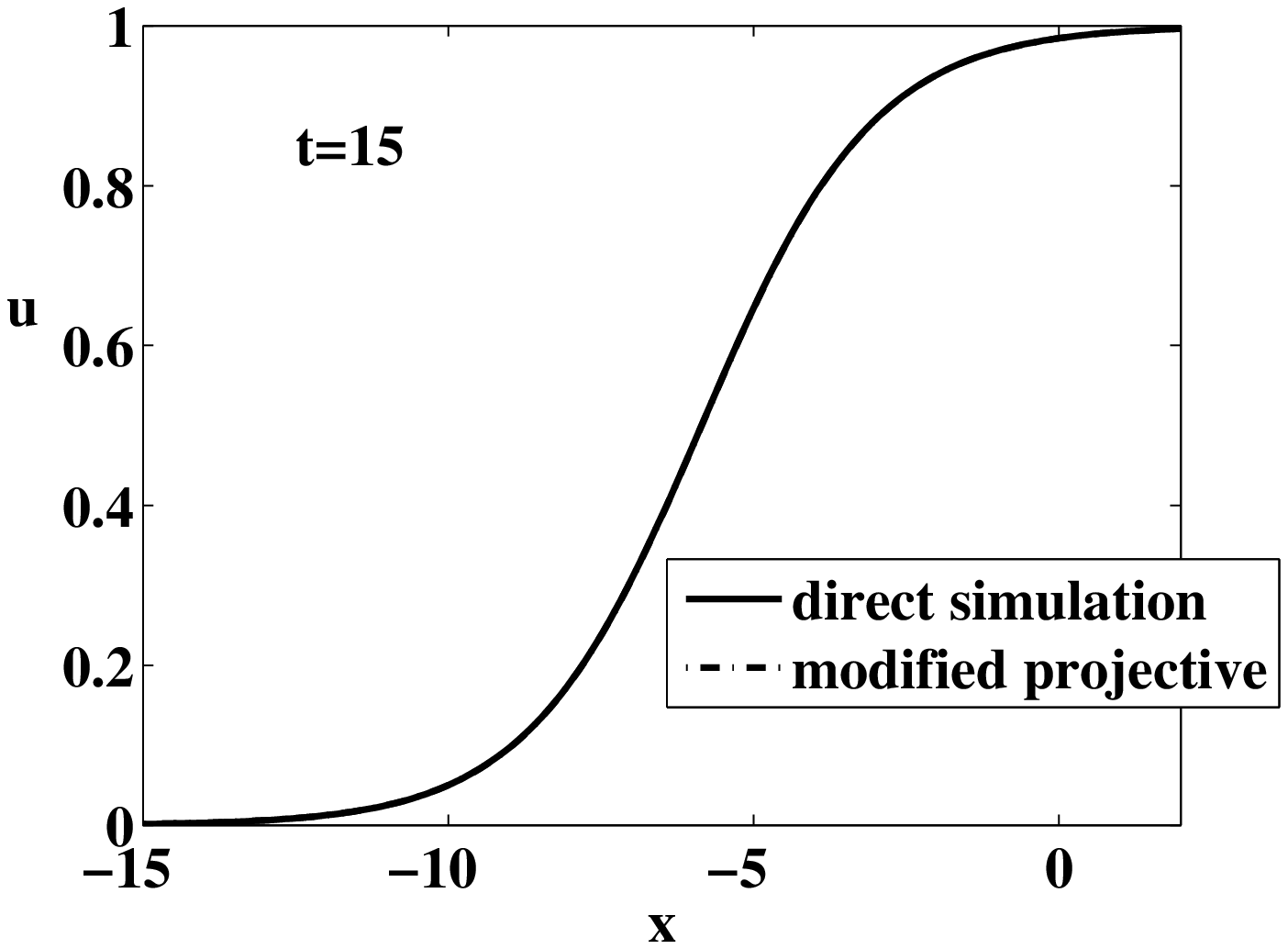} \\
(a)&(b)
\end{tabular}
\end{center}
\caption{Nagumo example: Solution obtained from direct simulation and (a) non co-traveling,
(b) co-traveling projective integration at time $t=15$ (see text).}
\label{ucompt15}
\end{figure}

A typical projective integration step requires the solutions $u_{1},u_{2}$ at two distinct reporting times
$t_{1},t_{2}$. 
To obtain the solution at projection time $t_{project}$,
we simply apply the Taylor expansion (\ref{Taylor_expansion}). 
We choose
two reporting times $t_{1}=0.1$ and $t_{2}=0.2$ which correspond to $2 \times 10^{3}$ steps of the
inner integrator. 
A projection time $t_{project}=0.5$ thus saves $3 \times 10^{3}$ inner integration steps.
If we use a projective method which ignores translational symmetry, the results are
manifestly inaccurate.
Fig.\ref{ucompt15}(a) compares the results of projective integration to those
of full direct simulation; taking translational invariance into account (see Fig.\ref{ucompt15}(b))
clearly shows the improved accuracy.

The template condition used to obtain these results  was:
\begin{equation}
\label{templ_FHN}
\int_{-30+c}^{30+c} \hat{u}(y) \D y \equiv \int_{-30}^{30} u(x+c,t) \D x = \int_{-30}^{30} u(x,0) \D x,
\end{equation}
implying that the integral of the {\em shifted} solution remains constant and equal to
the integral of the initial condition in the domain of interest.
Application of (\ref{templ_FHN}) at each reporting time $t_{1}$, $t_{2}$
produces the ``shifted'' solutions $\hat{u}_{1}$, $\hat{u}_{2}$ and the
corresponding shifts $c_{1}$, $c_{2}$.
The projection of $\hat{u}$ is obtained from
\begin{equation}
\hat{u}(t_{project})=
\hat{u}_{2}+(t_{project}-t_{2})
\frac{\partial{ \hat{u} }}{\partial{t}} \Bigg{|}_{t_{2}}
\approx
\hat{u}_{2}+(t_{project}-t_{2})
\frac{\hat{u}_{2}-\hat{u}_{1}}{t_{2}-t_{1}} 
\end{equation}
and the projection {\em of the shift} $c$ from
\begin{equation}
c(t_{project})=c_{2}+(t_{project}-t_{2})
\frac{\D c}{\D t} \Bigg{|}_{t_{2}}
\approx
c_{2}+(t_{project}-t_{2})
\frac{c_{2}-c_{1}}{t_{2}-t_{1}}.
\end{equation}
Finally, the full projected solution, reconstructed in physical space $x$ is recovered,
if desired, applying the inverse shift operator:
\begin{equation}
u(x,t_{project})=\hat{u}(x-c(t_{project}),t_{project}).
\end{equation}

\begin{figure}
\begin{center}
\includegraphics[width=0.6\linewidth]{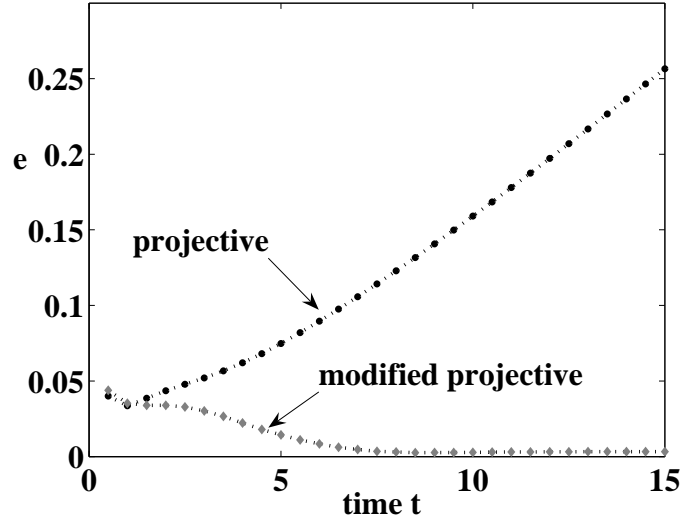}
\end{center}
\caption{Error evolution for the non co-traveling and for the modified
co-traveling projective integration in Fig.\ref{ucompt15}.}
\label{L2com}
\end{figure}

\begin{figure}
\begin{center}
\begin{tabular}{cc}
\includegraphics[width=0.49\linewidth]{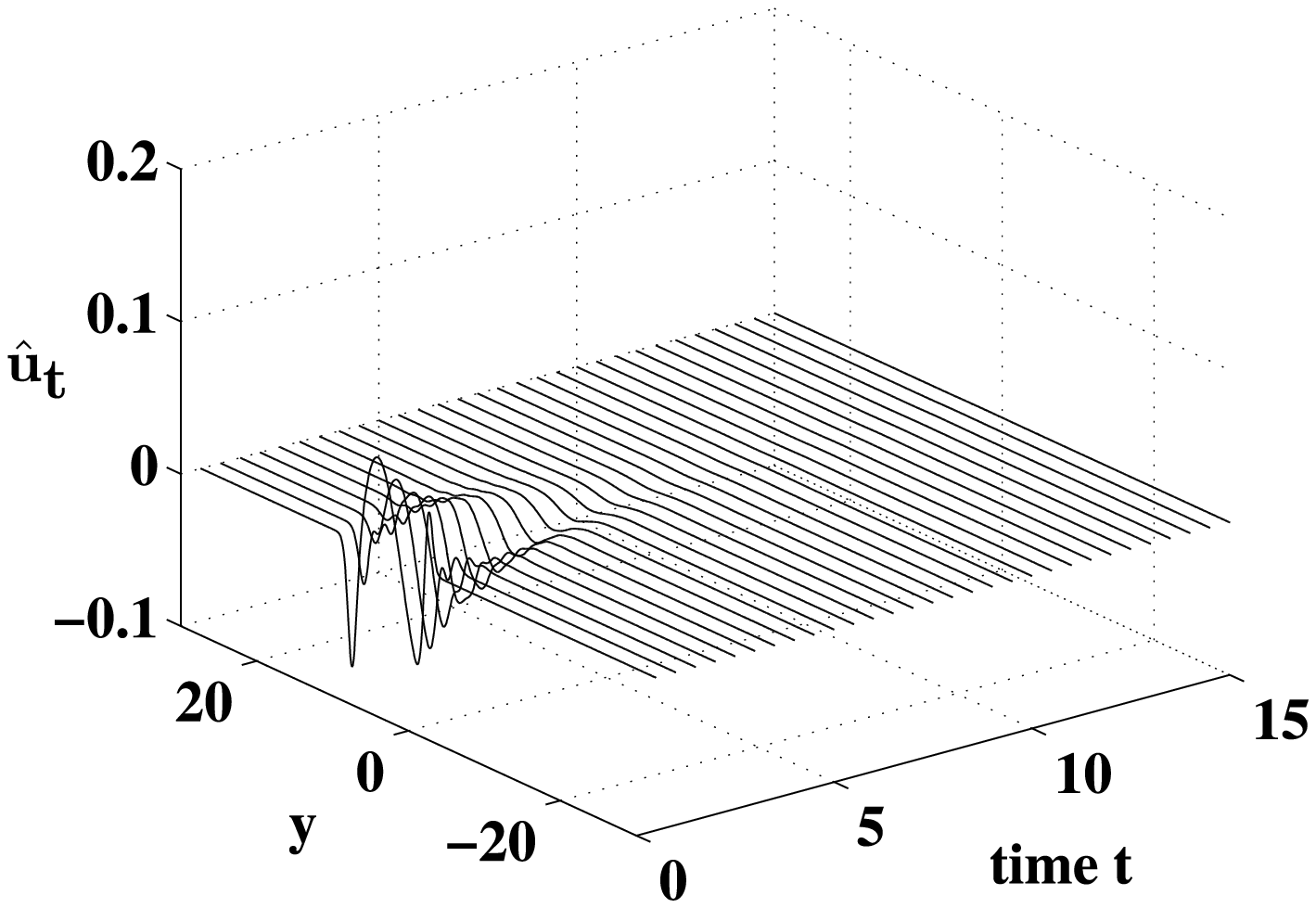} &
\includegraphics[width=0.49\linewidth]{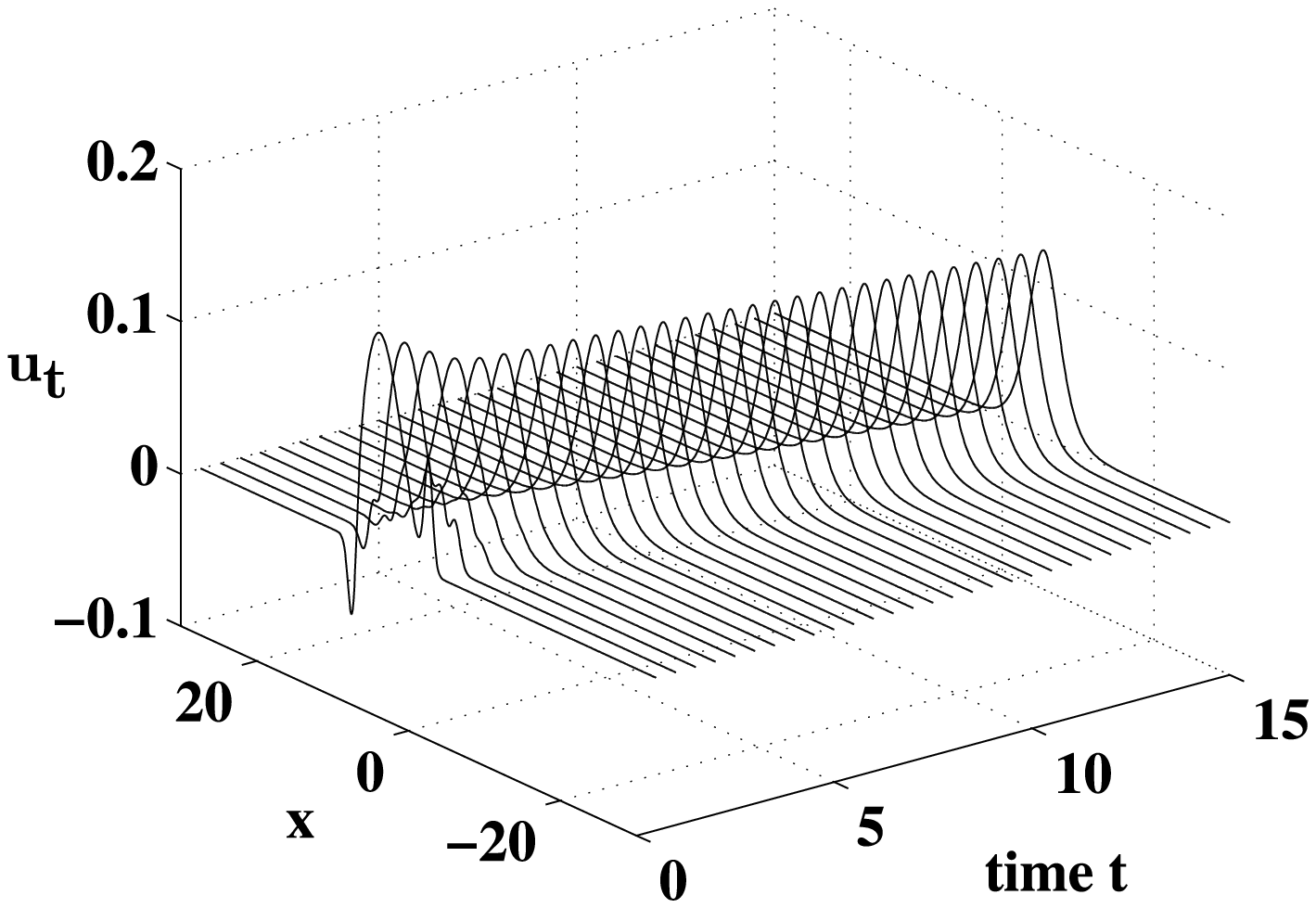} \\
(a)&(b)
\end{tabular}
\end{center}
\caption{Time derivative of the computed Nagumo solution (a) $\hat{u}$ and (b) $u$ up
to $t=15$.
The relatively high value of $u_{t}$ is the primary cause
for low accuracy of the obtained results when projective integration is applied
in a stationary frame.}
\label{dudt_comp}
\end{figure}

\begin{figure}
\begin{center}
\begin{tabular}{cc}
\includegraphics[width=0.45\linewidth]{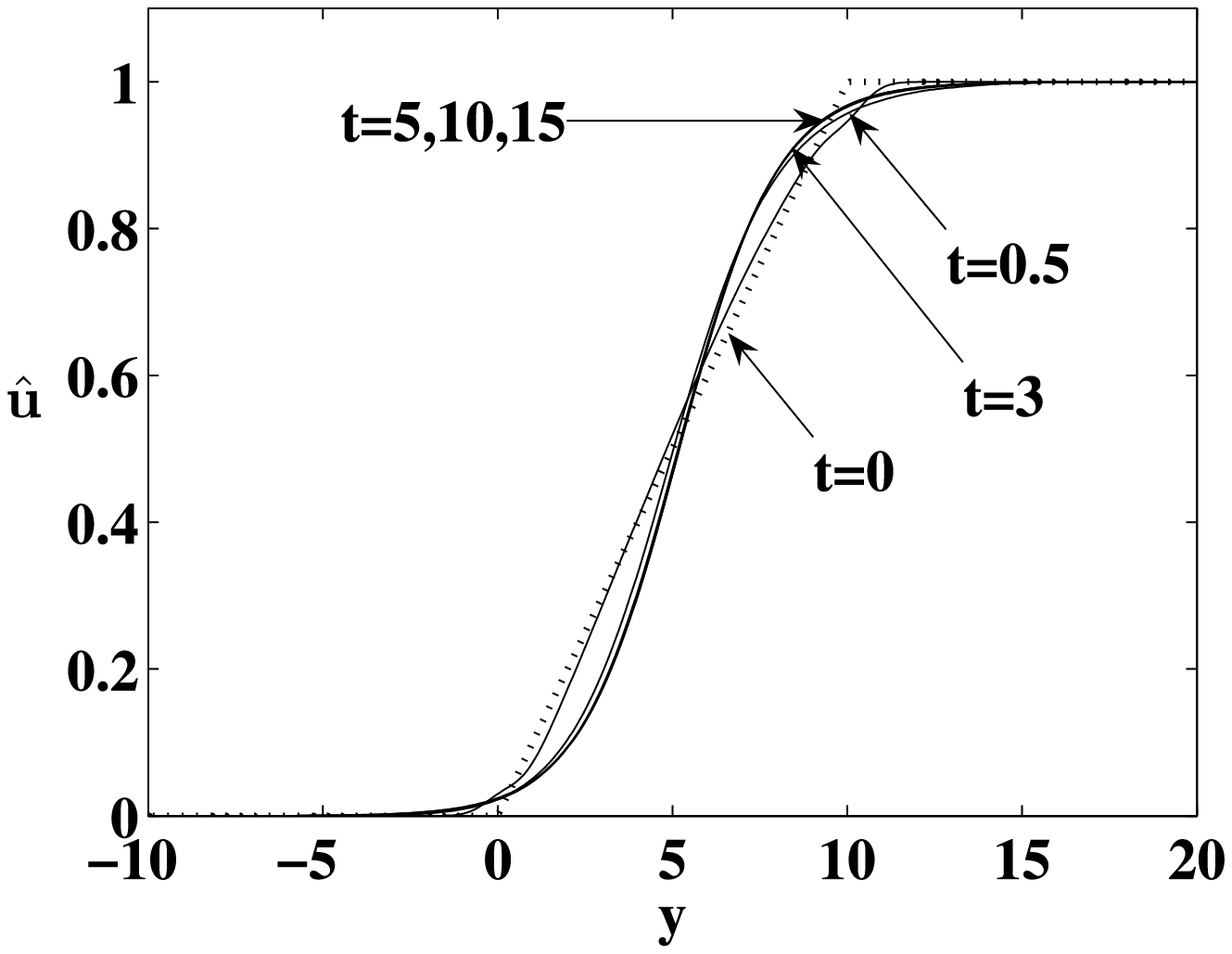}&
\includegraphics[width=0.45\linewidth]{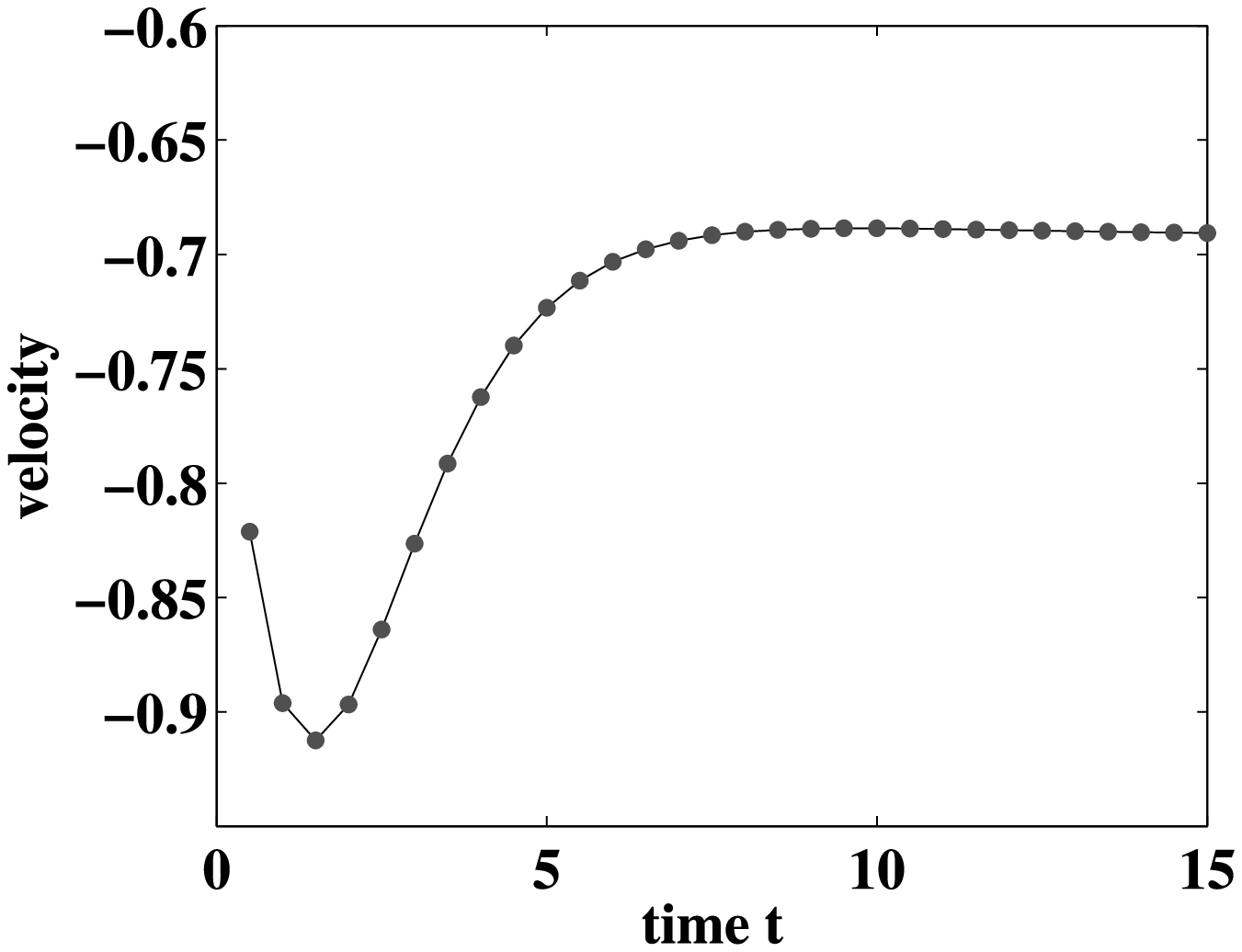} \\
(a)&(b)
\end{tabular}
\end{center}
\caption{(a) The shifted solution $\hat{u}$ evolves towards
a steady shape, the traveling wave of the Nagumo equation.
(b) The velocity converges to a constant value, $\D c / \D t \approx -0.69$, which agrees
with the theoretical value of $-0.693$.}
\label{cotraveling1}
\end{figure}

The benefit of projecting in a co-traveling frame is more clearly depicted
in Fig.\ref{L2com}, where we plot the error evolution ($L^{2}$ norm of the difference between
the solution obtained from direct simulation, $u_{direct}$ and the solution computed
from projective integration (co-traveling or not), $u_{PI}$ at the same time, $t$), i.e.:
\begin{equation}
e(t)=\left( \int_{-30}^{+30}(u_{direct}(x,t)-u_{PI}(x,t))^{2} \D x\right)^{1/2}.
\end{equation}
The error $e(t)$ resulting from the non co-traveling projective method
increases with time, while the results of 
the application of the modified projective integration scheme appear highly accurate.
The increased accuracy of the modified method is due to the slow evolution
(compared to the faster evolution of the unshifted solution)
(see Fig.\ref{dudt_comp}).
For the  co-traveling computations the solution after some time appears stationary (Fig.\ref{cotraveling1}(a));
this is (an approximation of) the stable Nagumo traveling wave.
Its (constant) speed (Fig.\ref{cotraveling1}(b)) can be approximated by 
\begin{equation}
\frac{\D c}{\D t} \approx \frac{c_{2}-c_{1}}{t_{2} -t_{1}}.
\end{equation}

\subsection{Coarse projective integration in a co-traveling frame - Kinetic Monte
Carlo simulation of a reaction-diffusion system}

We now apply the same methodology to
traveling problems for which the model simulations are conducted 
at a microscopic (stochastic, particle) 
level.

\begin{figure}[h]
\begin{center}
\includegraphics[width=0.5\linewidth]{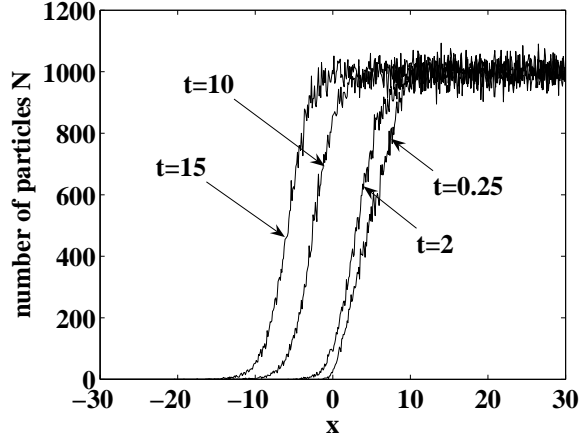}
\end{center}
\caption{Gillespie-based time evolution of a reaction-diffusion
 system motivated by Nagumo kinetics}
\label{Nagumo_straight}
\end{figure}

In our illustrative example the particle-based simulation is the Gillespie SSA
presented in Section \ref{FHNsim} applied to the same kinetic scheme. 
The one-dimensional domain of interest $[-30,30]$
is discretized  with
$J=601$ lattice sites. The distance, $h$, between two successive lattice sites
is equal to, $h=60/600=0.1$. 
The zero flux boundary conditions are 
incorporated applying a zero reaction rate for the 
``reactions'' (\ref{diffip1}) and (\ref{diffin1}) at 
sites $i=J$ and $i=1$ respectively, i.e. we do not allow the particles at $i=J$ to diffuse to the 
right and particles at $i=1$ to diffuse to the left.
At the deterministic limit the reaction rate constants are $k_{1}=1+\alpha=1.01$, $k_{-1}=1$ and
$k_{2}=\alpha=0.01$. The macroscopic diffusion coefficient, $D=1$,
corresponds to a diffusion rate constant $d=1/h^{2}$. In our computations,
we assume that the number of particles corresponding to dimensionless density,
$u=1$, is $N_{0}=1000$. The reaction parameters for the Gillespie code have been 
chosen consistently.
The number of particles at site $i$ is denoted by $N_{i}$, $i=1,\dots,601$.
The initial
condition is:
\begin{equation}
N_{i}=\left \{ \begin{array} {ccc}
0& \textrm{ for } & 1\le i \le 201, \\
100i& \textrm{ for } & 202 \le i \le 401, \\
1000& \textrm{ for } & 402 \le i \le 601.
\end{array}
\right.
\end{equation}
The results obtained from the kinetic Monte Carlo simulation are illustrated in Fig.\ref{Nagumo_straight},
where one can clearly see the formation of a (stochastic) traveling interface sweeping the one-dimensional
domain.
The stochastic simulation described above can be computationally intensive, especially if
one increases the number of particles. 
Such computations can be accelerated through
a {\em coarse} projective integration scheme;
translational invariance at the {\em coarse} (concentration field)
level should then be taken under consideration for more
accurate results.

We now describe the modified coarse projective integration scheme applied to
the Nagumo kinetics-motivated, kinetic Monte Carlo simulator.
The SSA is performed on $J=601$
lattice sites. In order to obtain a less noisy distribution we estimate a smoothened distribution 
in $n=101$ nodes, using the local averaging operator

\begin{equation}
M_{1}=\frac{1}{4} \sum_{k=1}^{k=4} N_{k}, \qquad
M_{i}=\frac{1}{7} \sum_{k=6i-9}^{k=6i-3} N_k, \textrm{ }i=2,\cdots,100, \qquad
M_{101}=\frac{1}{4} \sum_{k=598}^{k=601} N_{k}.
\end{equation}

We would like to approximate this distribution
in Fourier form; due to the boundary conditions, we
consider the difference distribution $f$ (in effect,
the spatial derivative) defined as:
\begin{equation}
f_{j}=\left \{ \begin{array} {ccc}
M_{j}-0& \textrm{ for } & j=1, \\
M_{j}-M_{j-1}& \textrm{ for } & j=2,..,n \\
\end{array}
\right.
\end{equation}
and its Fourier approximation \cite{Moeller:2003:EEDS}, i.e.:
\begin{equation}
\label{fourierexp}
f(x) \approx \frac{a_{0}}{2} + \sum_{k=1}^{k=K}\left(a_{k}\cos\left[k \frac{2\pi x}{L}\right]+b_{k}\sin\left[ k \frac{2\pi x}{L}\right] \right)
\end{equation}
where $L$ is the domain length ($L=60$).
The coarse variables in our computations are the first $K$ Fourier coefficients
of $f$. 
In the projective integration context we compute
these Fourier coefficients at two reporting times $t_{1},t_{2}$, approximate their time
derivatives at $t_{2}$, and extrapolate to the projection time $t_{project}$.
Such a computation, does not take into account the
translationally invariant character of the problem, leading to low accuracy results 
for relatively large steps. 
Application of the coarse projective scheme in a
co-traveling frame can capture the dynamics of the same system with enhanced accuracy,
even for relatively large projecting horizons, $T=t_{project}-t_{2}$.

We denote the shifted version of $f$, with $\hat{f}$, i.e.:
\begin{equation}
f(x)= \hat {f} (x-c).
\end{equation}
The Fourier coefficients $\hat{a}_{i},\hat{b}_{i}$ of $\hat{f}$ are then given
by:
\begin{equation}
\hat{a}_{i}=a_{i}\cos \left[i \frac{2\pi c}{L}\right]+b_{i}\sin\left[i \frac{ 2\pi c}{L}\right], \quad
\hat{b}_{i}=-a_{i}\sin\left[i \frac{2\pi c}{L}\right]+b_{i}\cos\left[i \frac{2\pi c}{L}\right].
\end{equation}
We apply the template condition 
\begin{equation}
\frac{\D}{\D c}\int_{-30}^{30}\hat{f}(x) \widetilde{T}(x) \, \D x = 0 \quad \Longleftrightarrow \quad
\frac{\D}{\D c}\int_{-30}^{30}{f}(x+c) \widetilde{T}(x) \, \D x = 0
\end{equation}
seeking maximum overlap between the shifted solution $\hat{f}$ and a template function $\widetilde{T}$.
Choosing the trigonometric template function $\widetilde{T}(x)=1-\cos[{2 \pi x}/{L}]$
reduces the template condition to:
\begin{equation}
\label{templatefourier}
\frac{\D \hat{a}_{1}}{\D c} = 0 \quad \Longrightarrow \quad
\frac{2 \pi c}{L}=\arctan\left[ \frac{b_{1}}{a_{1}} \right].
\end{equation}
The computation of the shifts $c_{1},c_{2}$ at reporting times $t_{1},t_{2}$, enables
the determination of the ``shifted'' Fourier coefficients $\hat{a}_{i},\hat{b}_{i}$ and their
projection at time $t_{project}$:
\begin{equation}
\label{aproject}
\hat{a}_{i}(t_{project}) = \hat{a}_{i}(t_{2})+
(t_{project}-t_{2})
\frac{\hat{a}_{i}(t_{2})-\hat{a}_{i}(t_{1})}{t_{2}-t_{1}} \quad
\textrm{for} \quad i=0,...,K
\end{equation}

\begin{equation}
\label{bproject}
\hat{b}_{i}(t_{project}) = \hat{b}_{i}(t_{2})+
(t_{project}-t_{2})
\frac{\hat{b}_{i}(t_{2})-\hat{b}_{i}(t_{1})}{t_{2}-t_{1}}
\quad \textrm{for} \quad i=1,...,K.
\end{equation}
The projected $f$ is then recovered by applying (\ref{fourierexp}).
The distribution $M$ of particles at the
$n$ nodes is then computed and the lifting procedure concludes by interpolating $M$ to
the $J$ lattice sites of the kinetic MonteCarlo time simulator
through MATLAB's intrinsic function {\tt interpft}. When this interpolation does not give 
an integer number of particles we round off.
This gives us the projected particle distribution in the co-evolving frame.
The spatial position of this distribution is determined by the projection
of the shift, $c_{project}$, according to
\begin{equation}
c_{project}=c_{2}+(t_{project}-t_{2})
\frac{c_{2}-c_{1}}{t_{2}-t_{1}}.
\end{equation}
Note that the projected solution does not necessarily exactly satisfy the template condition;
we circumvent this issue by calculating the  $\hat{b}_{1}(t_{project})$ 
from (\ref{templatefourier}) with $c=c_{project}$
and $a_{1}=a_{1}(t_{project})$.

The Fourier coefficients of the full, un-shifted solution, $f$  are
\begin{equation}
a_{i}(t_{project})=\hat{a}_{i}(t_{project})\cos\left[i \frac{2 \pi c_{project}}{L}\right]-
\hat{b}_{i}(t_{project})\sin\left[i \frac{2 \pi c_{project}}{L}\right]
\end{equation}
\begin{equation}
b_{i}(t_{project})=\hat{a}_{i}(t_{project})\sin\left[i \frac{2 \pi c_{project}}{L}\right]+
\hat{b}_{i}(t_{project})\cos\left[i \frac{2 \pi c_{project}}{L}\right],
\end{equation}
from which the projected particle distribution at the $j$ lattice sites 
can be obtained.

\begin{figure}
\begin{center}
\begin{tabular}{cc}
\includegraphics[width=0.5\linewidth]{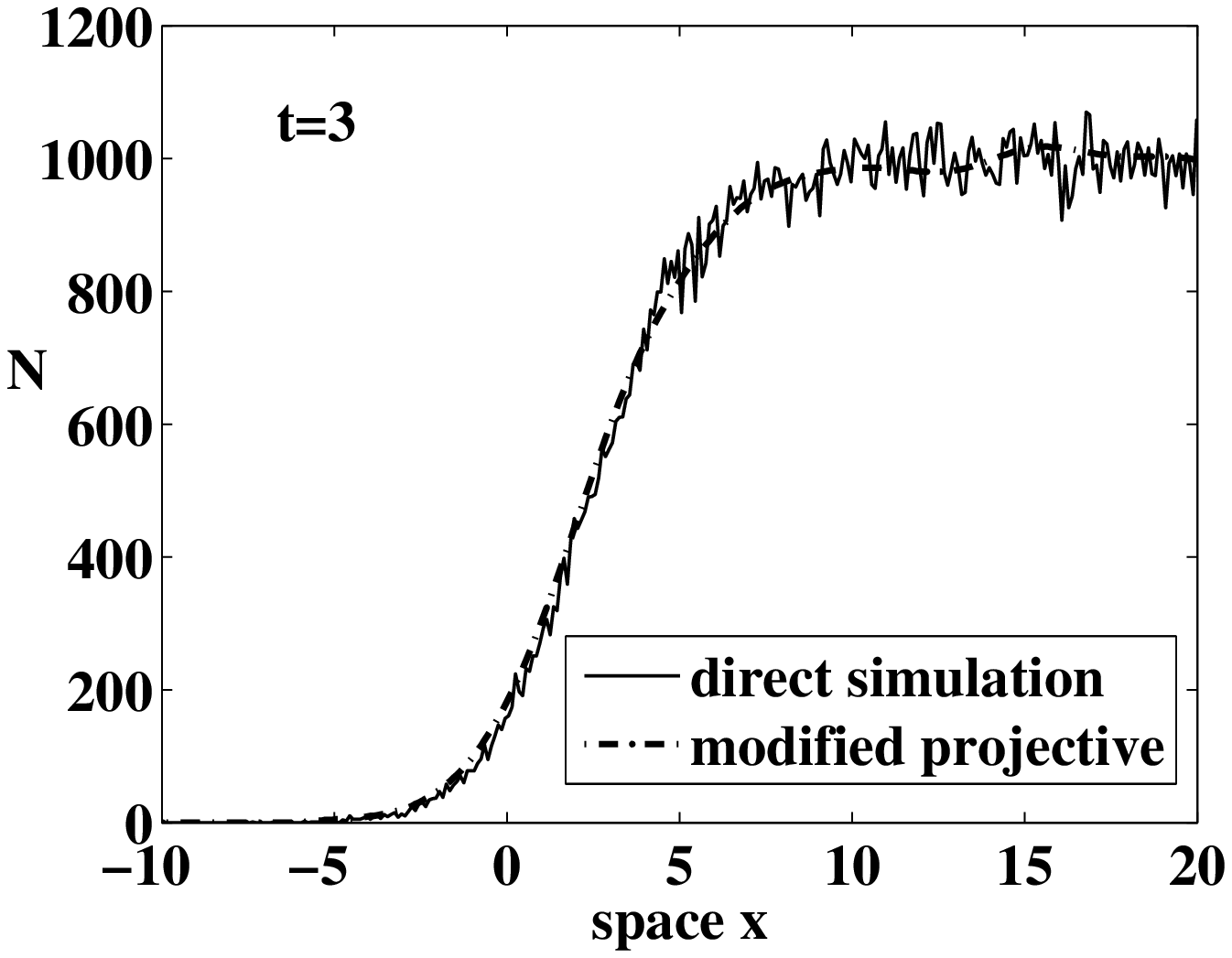} &
\includegraphics[width=0.5\linewidth]{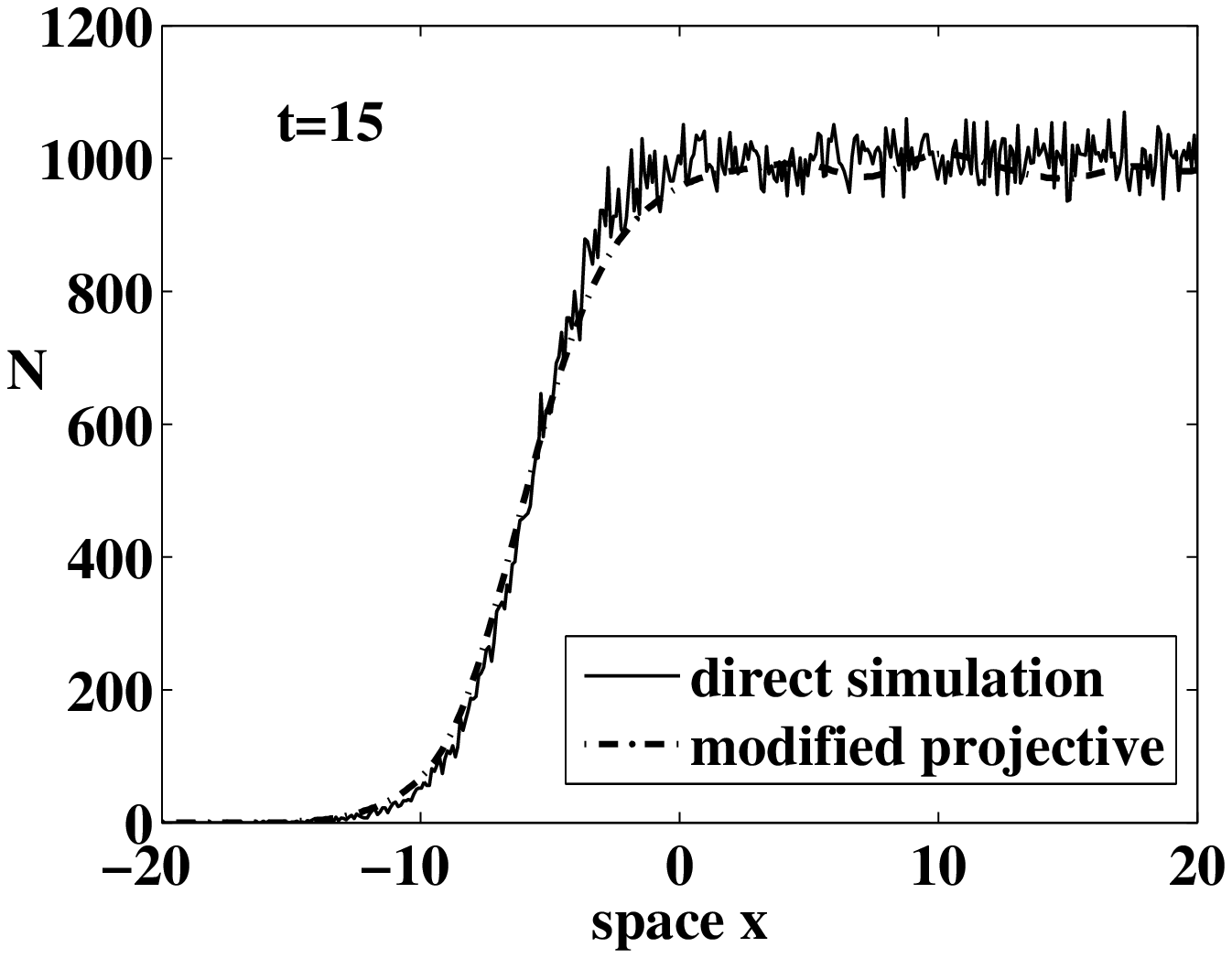} \\
(a)&(b)
\end{tabular}
\end{center}
\caption{Nagumo problem. Particle distribution obtained from direct simulation and coarse projective integration 
(a) after $t=3$ and (b) after $t=15$.}
\label{FHNGil}
\end{figure}

\begin{figure}
\begin{center}
\begin{tabular}{cc}
\includegraphics[width=0.5\linewidth]{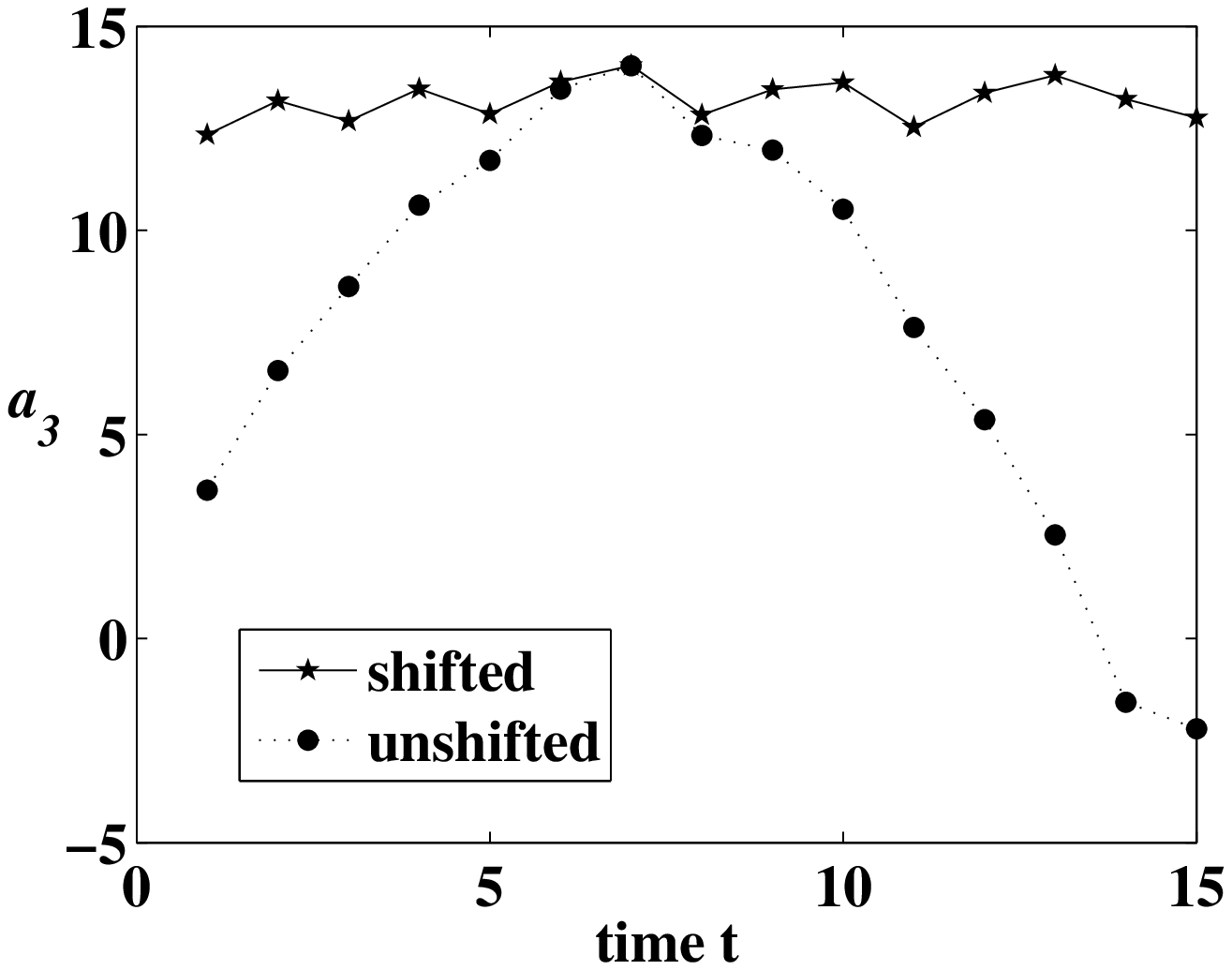} &
\includegraphics[width=0.5\linewidth]{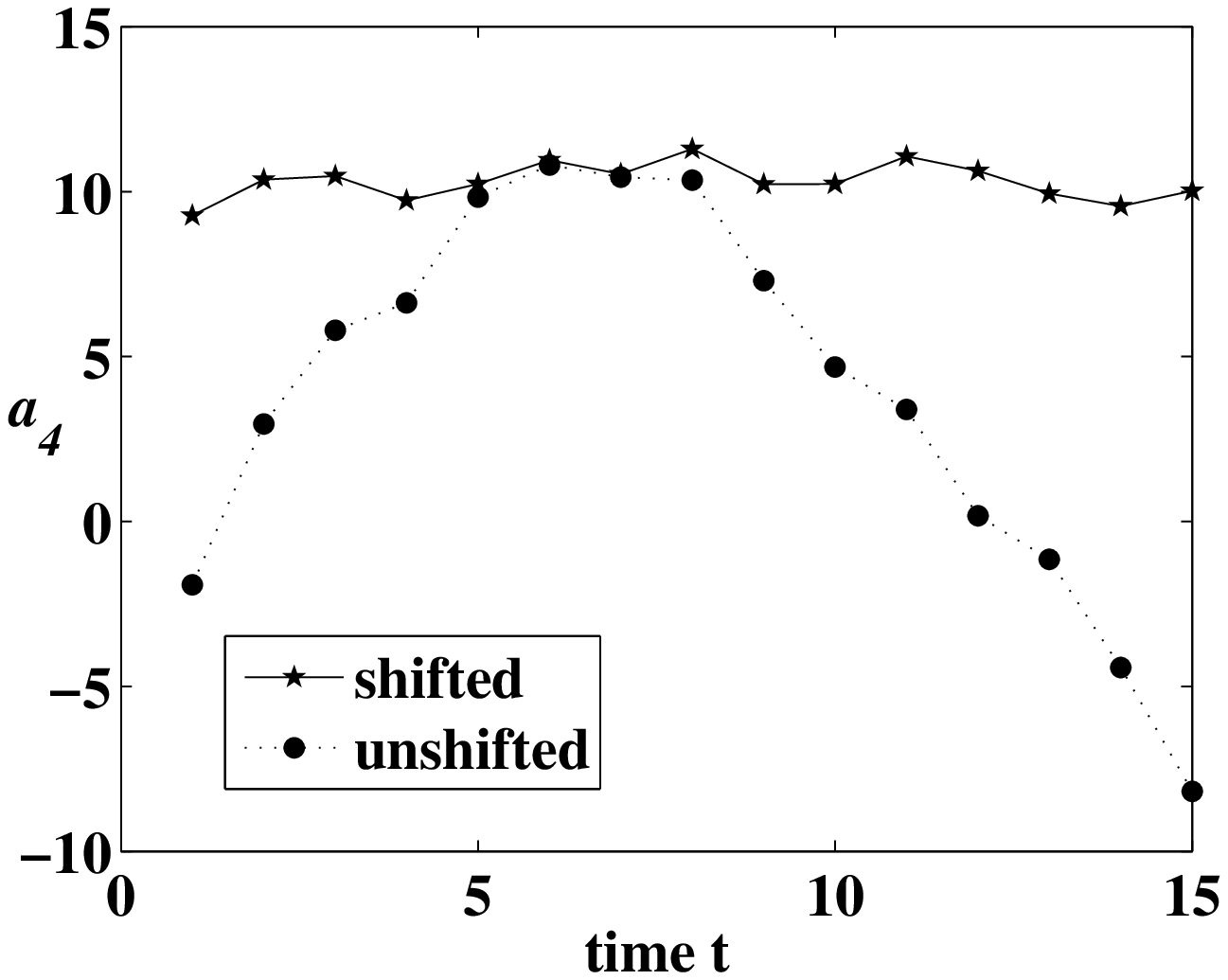} \\
(a)&(b)
\end{tabular}
\end{center}
\caption{Nagumo problem. Evolution of coarse variables ((a) 3rd Fourier coefficient $a_{3}$ and
(b) 4th Fourier coefficient $a_{4}$)
in the co-traveling frame (stars) and in a constant frame (dots).}
\label{FHNGilcoarse}
\end{figure}

Results of this template-based projective scheme are presented in Fig.\ref{FHNGil},
accurately capturing the (coarse) dynamics of the kinetic Monte Carlo Nagumo simulator, 
even for relatively large projection steps. 
The reporting times at each projective step were taken so that $t_{2}-t_{1}=0.25$,
the projection horizon was $t_{project}-t_{2}=0.5$ and the number of Fourier coefficients was $K=15$.
The coarse variables of the problem (Fourier coefficients) become essentially constant
in the co-traveling frame (see Fig.\ref{FHNGilcoarse}).


\section{Projective and Coarse projective integration in a dynamically rescaled frame: Diffusion}

\label{diffusionSIM}

In this section we study a projective scheme modified for scale
invariant systems, in particular systems that possess self-similar solutions.
Our illustrative example is simple one-dimensional diffusion,
both as a deterministic PDE and via a Monte Carlo-based simulation.

\subsection{One-dimensional diffusion - PDE example}
We study the simple mass diffusion equation
\begin{equation}
\label{diffusion}
u_{t}=u_{xx}.
\end{equation}
It possesses well-known self-similar solutions; we will exploit this
property, in order to perform relatively large projective steps accurately.
One can easily verify the
scale invariant character of (\ref{diffusion})
\begin{equation}
\label{scaleinvar_dif}
\operL_{x} \left(Bu\left(\frac{x}{A}\right)\right)=BA^{-2}\operL_{y} \left( u(y) \right)
\quad \mbox{where} \quad y=\frac{x}{A}.
\end{equation}

For our numerical computations we discretize the one-dimensional domain $x\in[-10,10]$
in $1001$ equidistant nodes (i.e., $\D x=0.02$). The dynamically renormalized
diffusion equation (along the lines of (\ref{coevolve_fin})) is
\begin{equation}
\frac{\partial{\hat{u}}}{\partial{\tau}}=
\operL_{y} \left(\hat{u}(y) \right) -\frac{1}{B}\frac{\D {B}}{\D {\tau}} \hat{u}
+\frac{1}{A}\frac{\D {A}}{\D {\tau}}y
\frac{\partial{\hat{u}}}{\partial{y}}.
\end{equation}
The spatial derivatives are approximated
with central finite differences and we consider zero flux boundary conditions.
The selected initial condition is
\begin{equation}
u(x,0)=\left \{ \begin{array} {ccc}
0& \textrm{ for } & |x|>1, \\
1& \textrm{ for } & |x|\leq1.
\end{array}
\right.
\end{equation}
At each step of the projective integration scheme, we choose two reporting
times $t_{1}$, $t_{2}$ and the solutions there, $u_{1}$, $u_{2}$
respectively. 

\begin{figure}
 \begin{center}
 \begin{tabular} {cc}
 \includegraphics[width=0.49\linewidth]{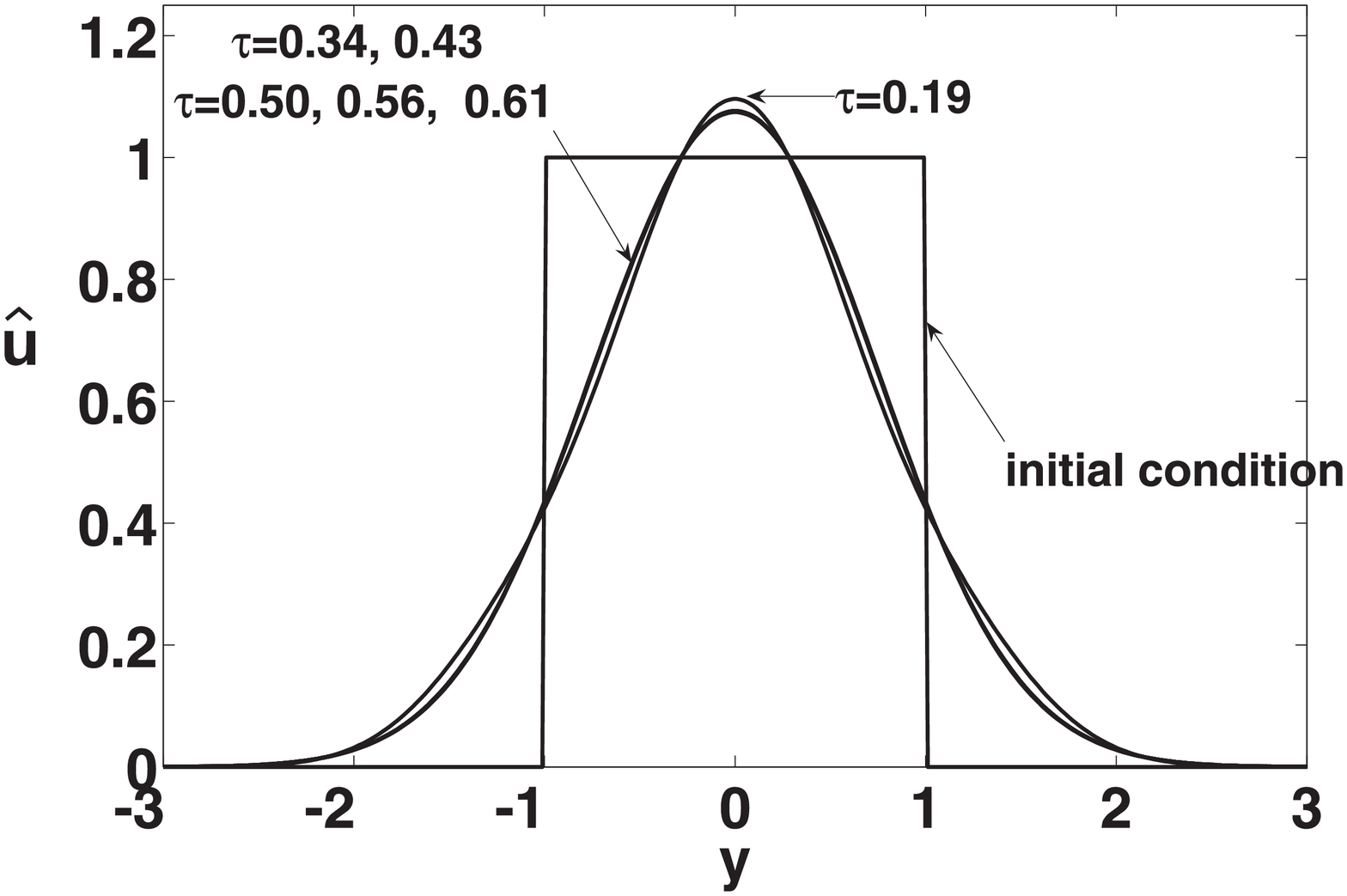} &
\includegraphics[width=0.49\linewidth]{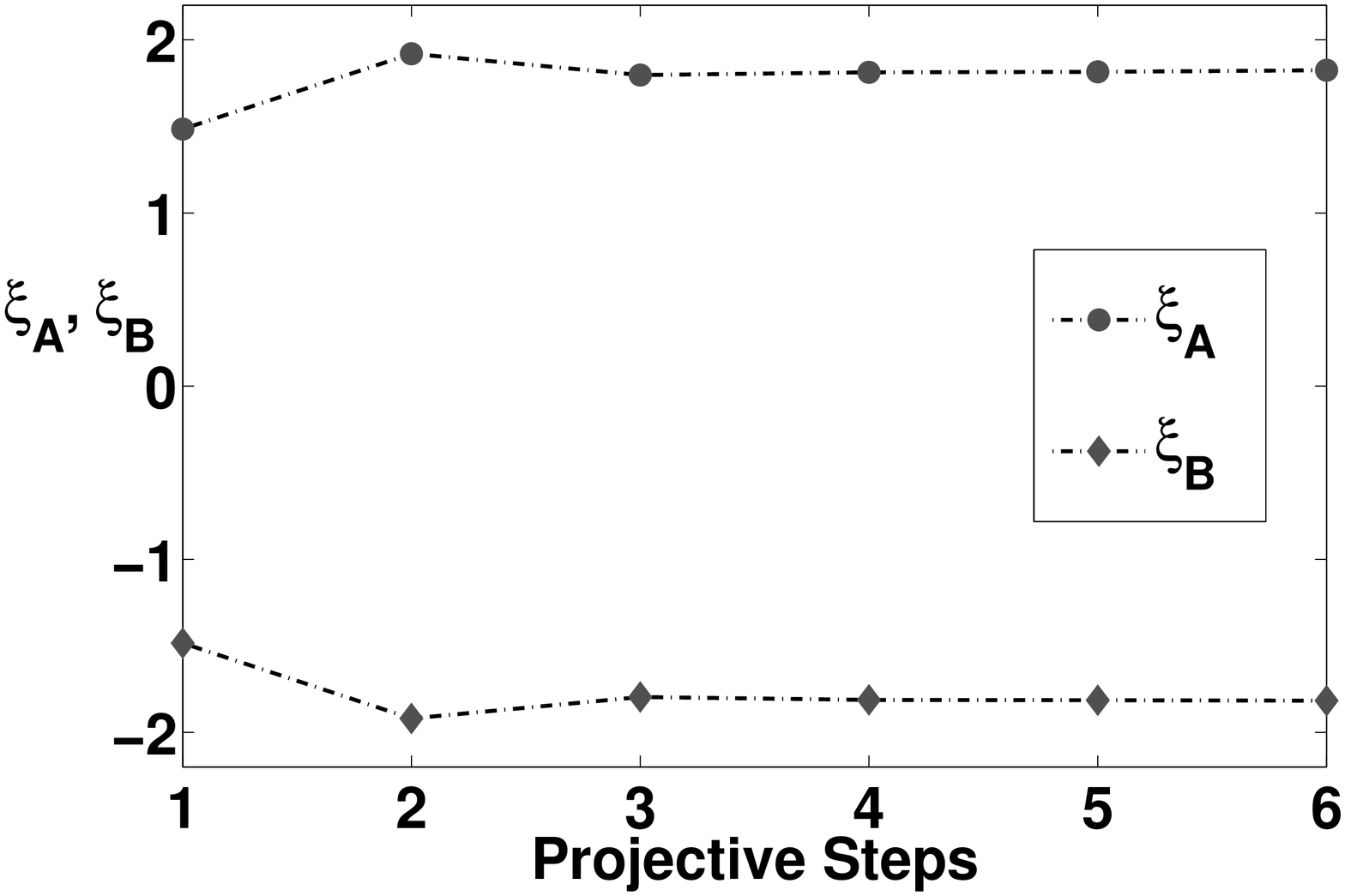} \\
  (a) & (b)
  \end{tabular}
  \end{center}
 \caption{Projective integration in a co-evolving (co-collapsing) frame applied to the 
 one-dimensional diffusion PDE example. (a) Instances of the evolution of rescaled
 solution $\hat{u}$ obtained at different -rescaled- projection times $\tau_{project}$.
The rescaled solution $\hat{u}$ converges to a steady state profile, corresponding
to a member of the self-similar family solutions. The initial condition
 is also depicted in the figure.
(b) Scale parameter $\xi_{A}=\frac{\D\log(A)}{\D\tau}$ and 
$\xi_{B}=\frac{\D\log(B)}{\D\tau}$ values computed at each projective 
step according to the procedure described in Section \ref{scale_invariance}.}

 \label{rescaleddiffusion}
 \end{figure}

A co-evolving frame formulation requires rescaling of the computed
solutions using
\begin{equation}
u(x,t)=B(\tau)\hat{u}\left( \frac{x}{A(\tau)},\tau(t)\right)
\end{equation}
as discussed in Section \ref{scale_invariance}. 
We can evaluate both the
rescaled solution, $\hat{u}$, and the scale factors $A,B$ at each reporting time step
by solving two template conditions. 
In this illustrative example the
template conditions chosen are
\begin{equation}
\int_{-\infty}^{+\infty}\hat{u}(y)\widetilde{T}_{1}(y) \D y=0 \quad \Longleftrightarrow \quad
\int_{-\infty}^{+\infty}\frac{1}{B}u(Ay)\widetilde{T}_{1}(y) \D y=0
\label{templA}
\end{equation}
for the template function
\begin{equation}
\widetilde{T}_{1}(y)=\left \{ \begin{array} {ccc}
-1& \textrm{ for } & |y|>1/2, \\
1& \textrm{ for } & |y|\leq 1/2
\end{array}
\right.
\end{equation}
and
\begin{equation}
\label{templB}
\int_{-\infty}^{+\infty}\hat{u}(y)\widetilde{T}_{2}(y) \D y=\mu \quad \Longleftrightarrow \quad 
\int_{-\infty}^{+\infty}\frac{1}{B}u(Ay)\widetilde{T}_{2}(y) \D y=\mu,
\end{equation}
where $\mu$ is a constant and the second template function is chosen as
$\widetilde{T}_{2}(y)=1$. 
The template condition (\ref{templB}) keeps 
the mass of the rescaled system constant, equal
to the initial mass.
Below we present results  
in a co-evolving frame, for $t_{2}-t_{1}=0.1$ and a projection step of
$t_{project}-t_{2}=0.2$ (thus economizing $10000$ time steps of the {\em inner} Euler integrator).
The evolution of $\hat{u}$ is depicted in Fig.\ref{rescaleddiffusion}; 
a stationary profile is approached  after some 
``rescaled'' time, $\tau$; this profile is a member
of the family of self-similar solutions of the diffusion equation. 
The scale parameters
$\xi_{A}=\frac{\D\log(A)}{\D\tau}$ and $\xi_{B}=\frac{\D\log(B)}{\D\tau}$
shown in the same figure
also approach stationarity. 

\begin{figure}
 \begin{center}

 \includegraphics[width=0.5\linewidth]{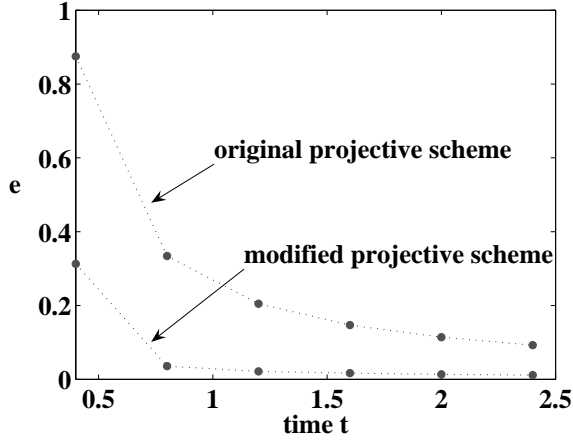}
  \end{center}
 \caption{Comparison of error evolution for the simulations in 
Fig.\ref{rescaleddiffusion} when the projective integration and the modified projective integration
is applied.
The modified projective integration
algorithm application manages to produce accurate results even at the early stages, where
the solution still evolves towards its self-similar shape.}
 \label{L2compdif}
 \end{figure}

Comparison of the errors 
\begin{equation}
e(t)=\int_{-10}^{10} \left(u_{direct}(x,t)-u_{PI}(x,t)\right)^{2} \D x
\end{equation}
for
projective integration in a co-evolving frame
with those for unmodified projective integration  (shown in Fig.\ref{L2compdif})
illustrates the advantage of projecting in a dynamically renormalized frame;
$u_{direct}$ is the solution obtained from direct simulation of (\ref{diffusion}). The
errors are computed for the reconstructed solutions $u_{PI}$. 

\begin{figure}
 \begin{center}
 \begin{tabular} {cc}
 \includegraphics[width=0.49\linewidth]{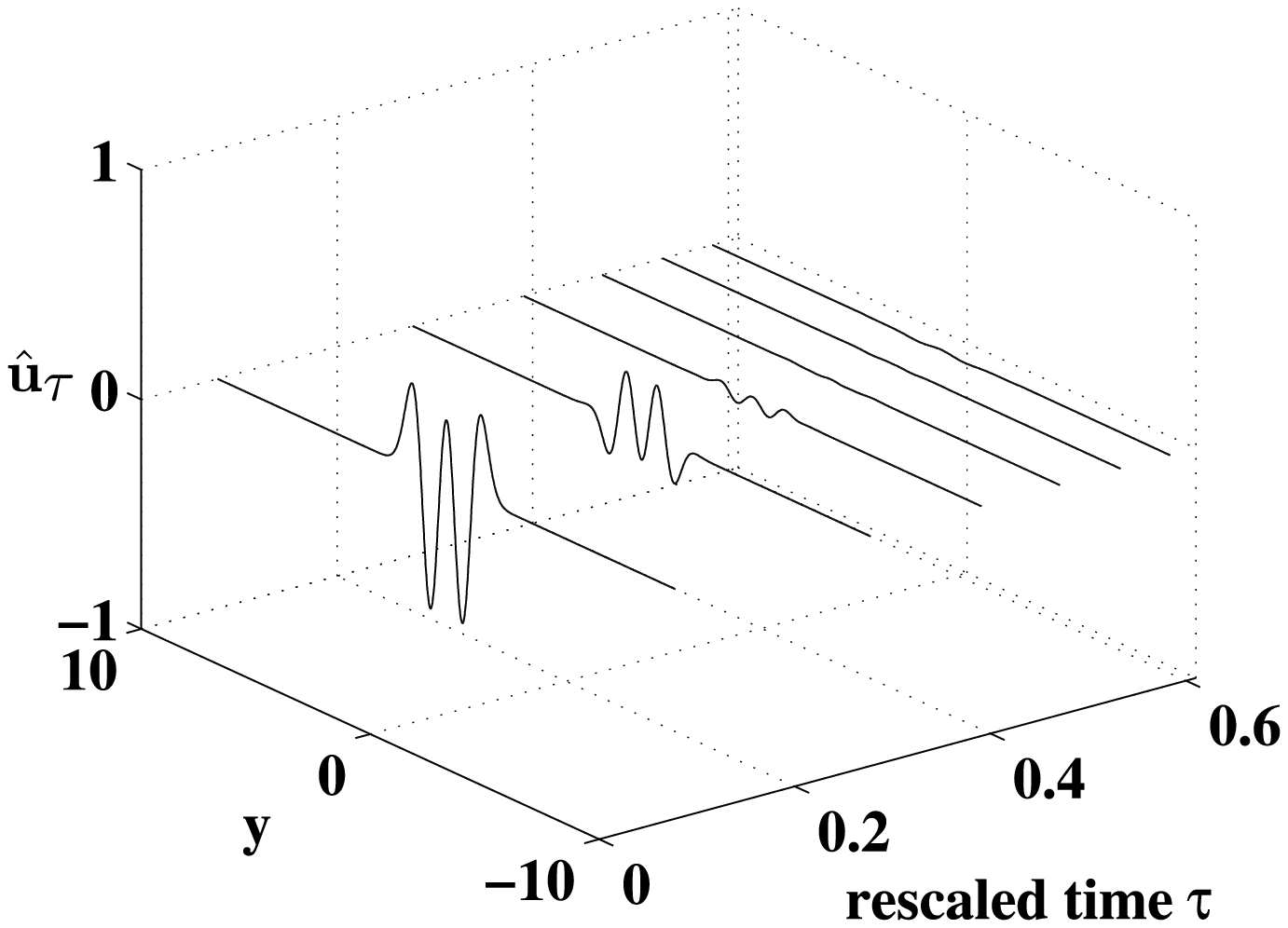} &
\includegraphics[width=0.49\linewidth]{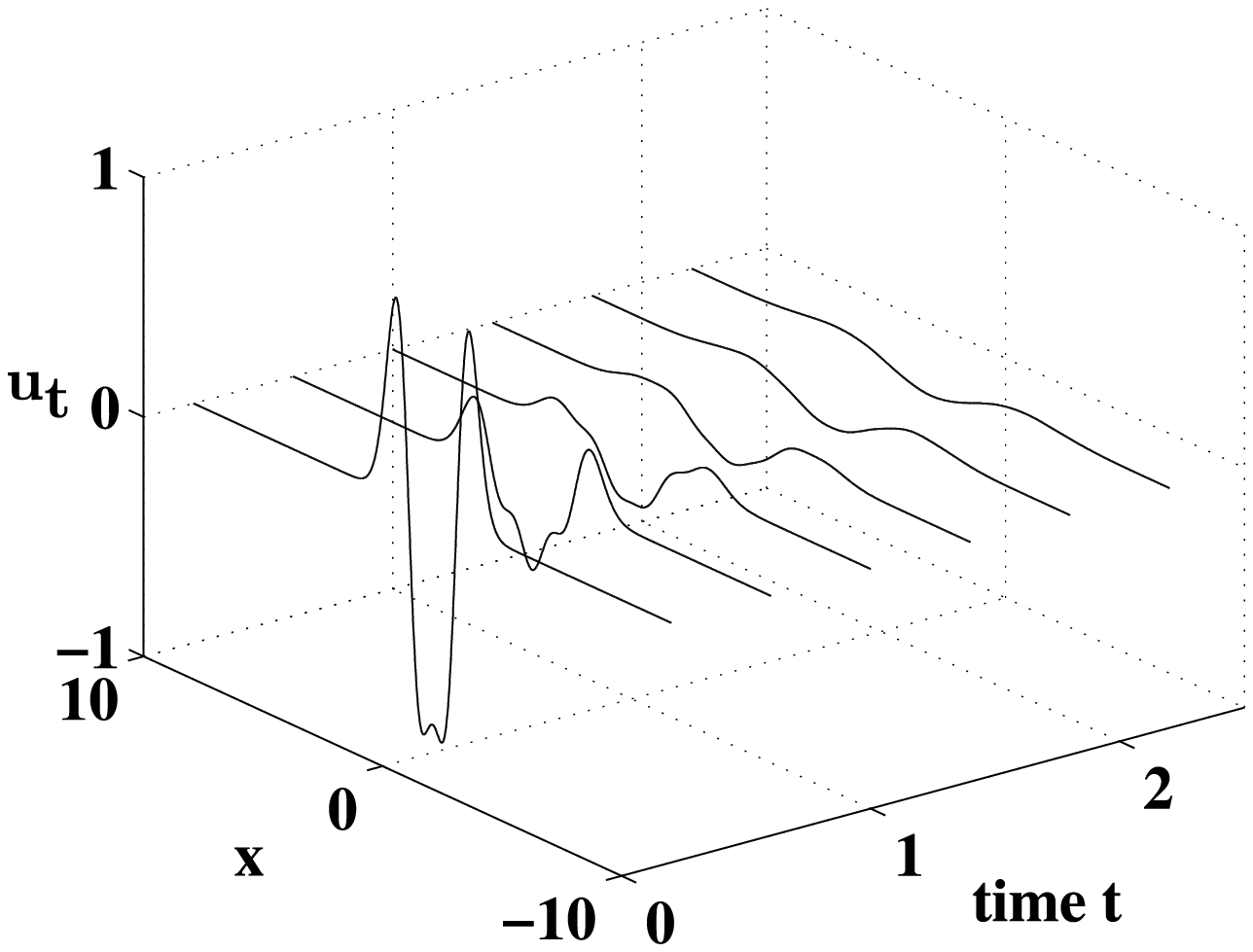} \\
  (a) & (b)
  \end{tabular}
  \end{center}
 \caption{One-dimensional diffusion equation. 
(a)``Time'' derivative of rescaled solution $\hat{u}$ at 
$\tau_{2}$-reporting times, using the (\ref{templA}) and
(\ref{templB}) template conditions. (b)
Time derivative of $u$ as obtained from the application of the original
projective integration scheme at $t_{2}$-reporting times. The relatively high values of
time derivative $u_{t}$ is the main reason for the failure of the method to
 capture the correct dynamics of (\ref{diffusion}).}
 \label{dudtcompdif}
 \end{figure}

Once more, the extra accuracy can be attributed to the slower evolution in the
dynamically renormalized frame (see Fig.\ref{dudtcompdif}).

\subsection{Random walker simulation of one-dimensional diffusion - Coarse Projective Integration}

In this section we present ``renormalized coarse projective integration"
applied to a Monte Carlo algorithm simulating diffusion in a
population of $10^{6}$ random walkers in one space dimension.

The macroscopic observable in this case is the cumulative
distribution function (CDF) of the particle positions denoted by $f$.
The domain of interest $[-10,10]$ is discretized into 1001
equally spaced nodes ($\D x=0.02$).
The CDF is then determined as a function of the discretized spatial domain
$\{x_{1},...,x_{1001} \}$. 
Each particle, $i$, is described by its position
$X_{i}$. 
The Monte Carlo time step is
$\delta t=0.0001$ and during each time step each particle will randomly 
move left or right with equal probability by an increment $\delta X=\sqrt{2\delta t}$.

\begin{figure}
\begin{center}
\includegraphics[width=0.7\linewidth]{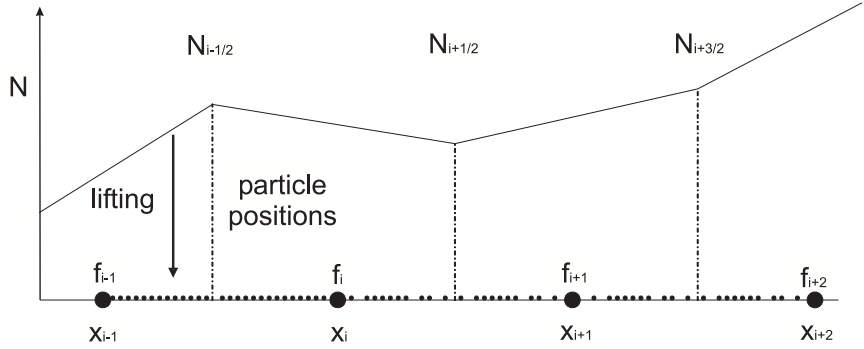}
\end{center}
\caption{The piecewise linear distribution of density $N$ corresponds to a microscopic realization of
particle positions. The distribution $N$ is evaluated as the spatial derivative of the CDF.}
\label{lifting_process}
\end{figure}

Denoting CDF at mesh point $x_{i}$ as $f_{i}$, 
the probability density function of particles
is evaluated by
\begin{equation} \label{N_density}
N_{i-1/2}=\frac{f_{i}-f_{i-1}}{x_{i}-x_{i-1}}
\end{equation}
where $N_{i-1/2}$ is the macroscopic density of particles at the midpoint 
$[x_{i-1}+ x_{i}]/2$. {\it Lifting} -- i.e. construction of  a microscopic state
consistent with density (\ref{N_density}) -- is done as follows.
The particles are placed in space so that their density piecewise
linearly interpolates the midpoint values (see Fig.\ref{lifting_process}).

We assume that an evolution equation for the CDF of particle positions $f(x)$ exists:
\begin{equation}
\label{cdf_evol_dif}
\frac{\partial{f}}{\partial{t}}=\operL_{x}(f).
\end{equation}
Before applying the coarse projective scheme in a co-evolving framework
we should test the scale invariance of the unknown differential operator
$\operL_{x}$ and extract its scaling exponents.
In our computations there is no amplitude scaling 
since $f$ is a CDF.
The operator should satisfy 
\begin{equation}
\label{scale_invariance_diffusion}
\operL_{x} \left( f\left(\frac{x}{A}\right)\right) = A^{a} \operL_{y} \left(f(y)\right)
\quad \mbox{where} \quad y=\frac{x}{A}
\end{equation}
for any $A$.
Despite the fact that the explicit formulation of $\operL_{x}$
is unknown we can estimate its action on test distribution $f(x)$ through
the computation of $\frac{\partial{f}}{\partial{t}}$.
We perform short
computational experiments to estimate the action of operator $\operL_{x}$
as follows:

\leftskip 1.1cm

\bigskip
\noindent
(1) We select a test function $\varphi_{0}$ and a positive constant $A$. 

\smallskip
\noindent
(2) We initialize the kinetic Monte Carlo 
simulator so that the CDF of particles is equal to $\varphi_{0}$ (using the
lifting procedure in Fig.\ref{lifting_process}).
We run the kinetic Monte Carlo simulator for a relatively short time interval
(in macroscopic terms) $\delta T$.
We obtain the new CDF $\varphi_1$
and estimate the time derivative of $\varphi$ from the expression
\begin{equation}
\label{time_der_f}
\frac{\partial{\varphi}}{\partial{t}} \approx \frac{\varphi_1-\varphi_{0}}{\delta T}.
\end{equation}

\smallskip

\noindent
(3) We initialize the kinetic Monte Carlo 
simulator so that the CDF of particles is equal to $\hat{\varphi}_{0}(x)=\varphi_{0}(x/A)$ (using the
lifting procedure in Fig.\ref{lifting_process}).
We run the kinetic Monte Carlo simulator for a short time interval $\delta T$.
We obtain the new CDF $\hat{\varphi}_1$
and estimate the time derivative of $\hat{\varphi}$ similarly as in (\ref{time_der_f}).

\smallskip
\noindent
(4)  We estimate the value of exponent $a$ in (\ref{scale_invariance_diffusion}) 
by minimizing the residual
\begin{equation}
\label{minimization}
R(a)=\left|\left| \frac{\partial{\hat{\varphi}}}{\partial{t}}\Big{|}_{Ax}-A^{a}
 \frac{\partial{\varphi}}{\partial{t}}\Big{|}_{x} \right|\right|^{2},
\end{equation}
where $||\cdot||$ denotes the standard Euclidean norm.

\bigskip
\leftskip 0cm

\begin{figure}
 \begin{center}
 \begin{tabular} {cc}
 \includegraphics[width=0.5\linewidth]{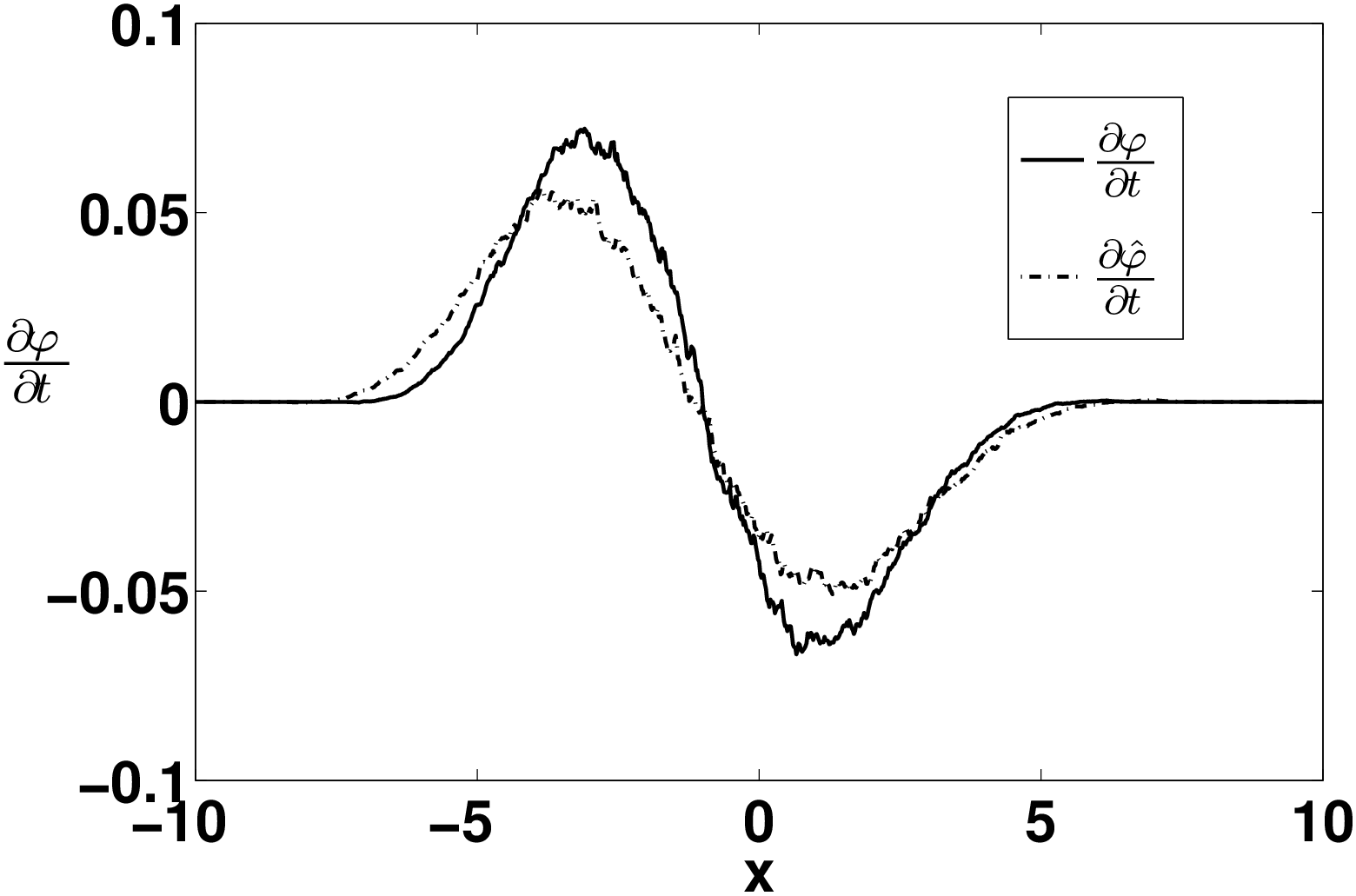} &
\includegraphics[width=0.5\linewidth]{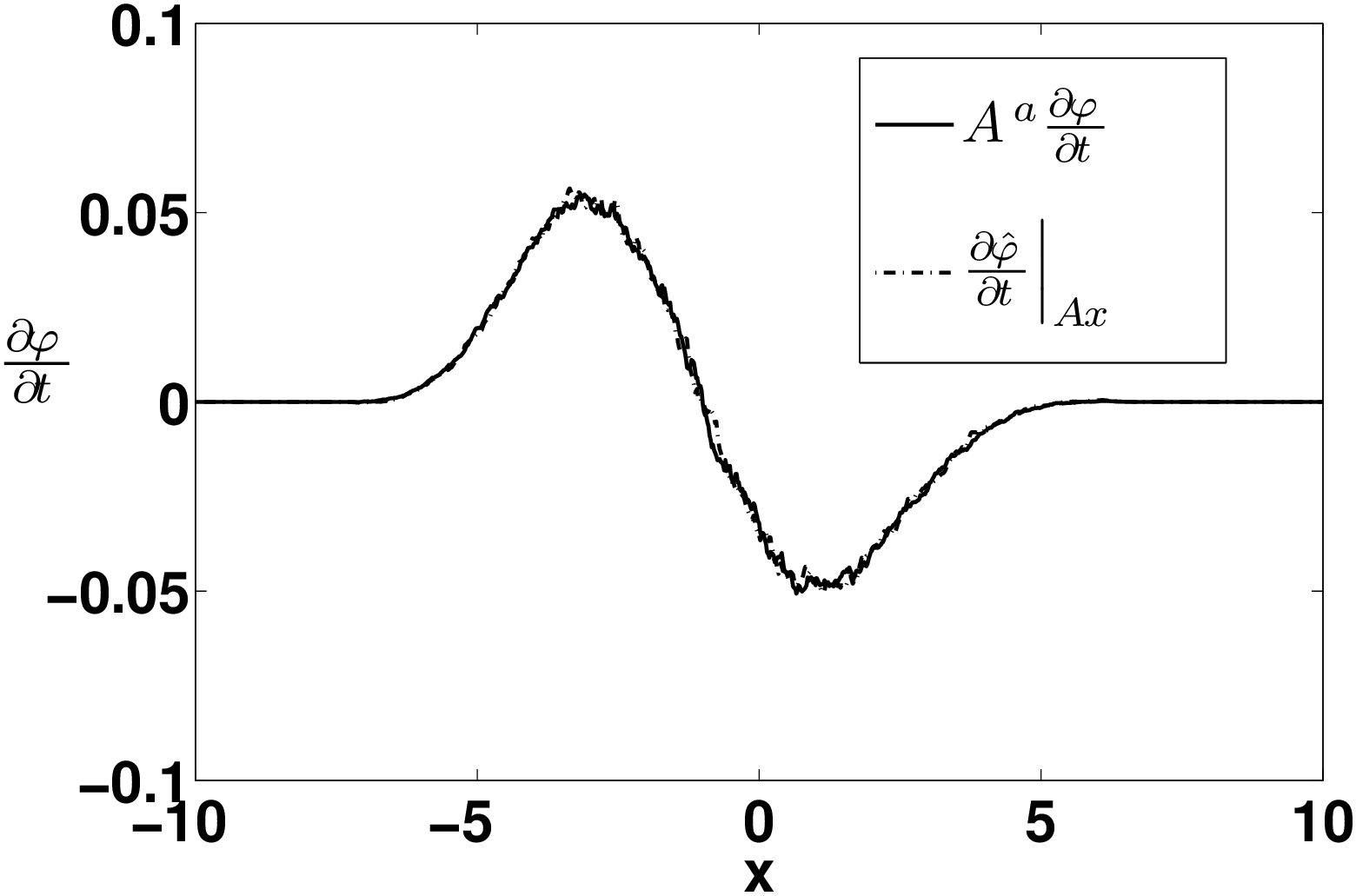} \\
 (a) & (b)\\
  
  \end{tabular}
  \end{center}
 \caption{(a) Time derivative of test function $\varphi_{0}$ given by (\ref{betadistr}) (solid line) and the rescaled
test function $\hat{\varphi}_{0}(x)=\varphi_{0}(\frac{x}{A})$ (dashed line); parameters of algorithm
(1) -- (4) are given in the text.
 (b) Testing the scale invariance property of operator $\operL_{x}$
 (see (\ref{scale_invariance_diffusion})) scale invariance property. The value of scale
 exponent is $a = -1.97 \approx -2$.}
 \label{verification_figure}
 \end{figure}

\noindent
In this case we
use:
\begin{equation}
\label{betadistr}
\varphi_{0,i}=\frac{1}{B(\gamma,\delta)} 
\int_{0}^{x_{i}/20+1/2} \zeta^{\gamma-1} (1-\zeta)^{\delta-1} \D \zeta,
\end{equation}
where $\varphi_{0,i}$ is the value of test function $\varphi_{0}$ at mesh point $x_{i}$,
$B(\gamma,\delta)$ is the Beta function with parameters $\gamma=8$ and 
$\delta=10$, $\delta {T}=0.01$ and $A=1.15$.
The results are shown in Fig.\ref{verification_figure}.
The procedure described above can be performed for different
values of $A$ and different test functions $\varphi_{0}$.
The value of the exponent $a$ minimizing the residual (\ref{minimization})
was in all tested cases close to $-2$. 
It confirms
the operator's $\operL_{x}$ scale invariance property
($\frac{\partial{\hat{\varphi}}}{\partial{t}}\Big{|}_{Ax}$ almost coincides with
$A^{a}\frac{\partial{\varphi}}{\partial{t}}\Big{|}_{x}$).
The value of the scaling exponent $a=-2$ is used in the proposed co-evolving
projective integration scheme.

In our modified coarse projective scheme we evaluate the 
macroscopic observables at $k=2$ distinct reporting times 
$t_{1}$, $t_{2}$ with corresponding CDFs $f_{1}$, $f_{2}$.
At each step of the projective integration
scheme the time step is $t_{2}-t_{1}=0.05$;
the projection step is $t_{project}-t_{2}=0.1$.
The initial CDF $f_{0}$ at mesh point $x_{i}$ is given by
\begin{equation}
f_{0}(x_{i})\equiv f_{0,i}=\left \{ \begin{array} {ccc}
0& \textrm{ for }&  1\le i \le 451, \\
{[i-451]}/100 & \textrm{ for } & 452\le i \le 551, \\
1 & \textrm{ for } & 552 \le i \le 1001.
\end{array}
\right.
\end{equation}
Both the scale factor $A$ and the rescaled solution $\hat{f}$, where
\begin{equation}
\hat{f}(x)=f(Ax),
\end{equation}
are obtained from the application of the template condition:
\begin{equation}
\label{templateA}
\hat{f}(\zeta_{2})-\hat{f}(\zeta_{1})=\nu \quad \Longrightarrow \quad
f(A\zeta_{2})-f(A\zeta_{1})=\nu,
\end{equation}
where $\zeta_{1}$, $\zeta_{2}$ are given real numbers and $\nu$ 
is constant.
This template condition, enforces a constant number of particles in
interval $[\zeta_{1},\zeta_{2}]$. Let us note that the CDF is defined only at mesh points 
$\{x_{1},...,x_{1001} \}$. Whenever template condition (\ref{templateA}) 
requires values of $f$ outside mesh points $\{x_{1},...,x_{1001} \}$, we use 
linear interpolation.
For our computations we chose $\zeta_{1}=-0.25$ and
$\zeta_{2}=0.5$, while the constant $\nu$ is evaluated from the initial condition $f_{0}$, i.e.
$f_{0}(\zeta_{2})-f_{0}(\zeta_{1})=\nu$.

As in the PDE diffusion example, we choose to evolve the rescaled CDF in rescaled
space $y$ and time $\tau$. 
The temporal Taylor expansion is performed in terms of $\tau$; we therefore
need to evaluate the $\tau_{1}$, $\tau_{2}$ and $\tau_{project}$ values
corresponding to $t_{1}$, $t_{2}$, $t_{project}$.
The procedure was reported in Section \ref{scale_invariance}; note that here
the $B$ scaling is omitted since $f$ is a CDF.

\begin{figure}
 \begin{center}
 \begin{tabular} {cc}
 \includegraphics[width=0.49\linewidth]{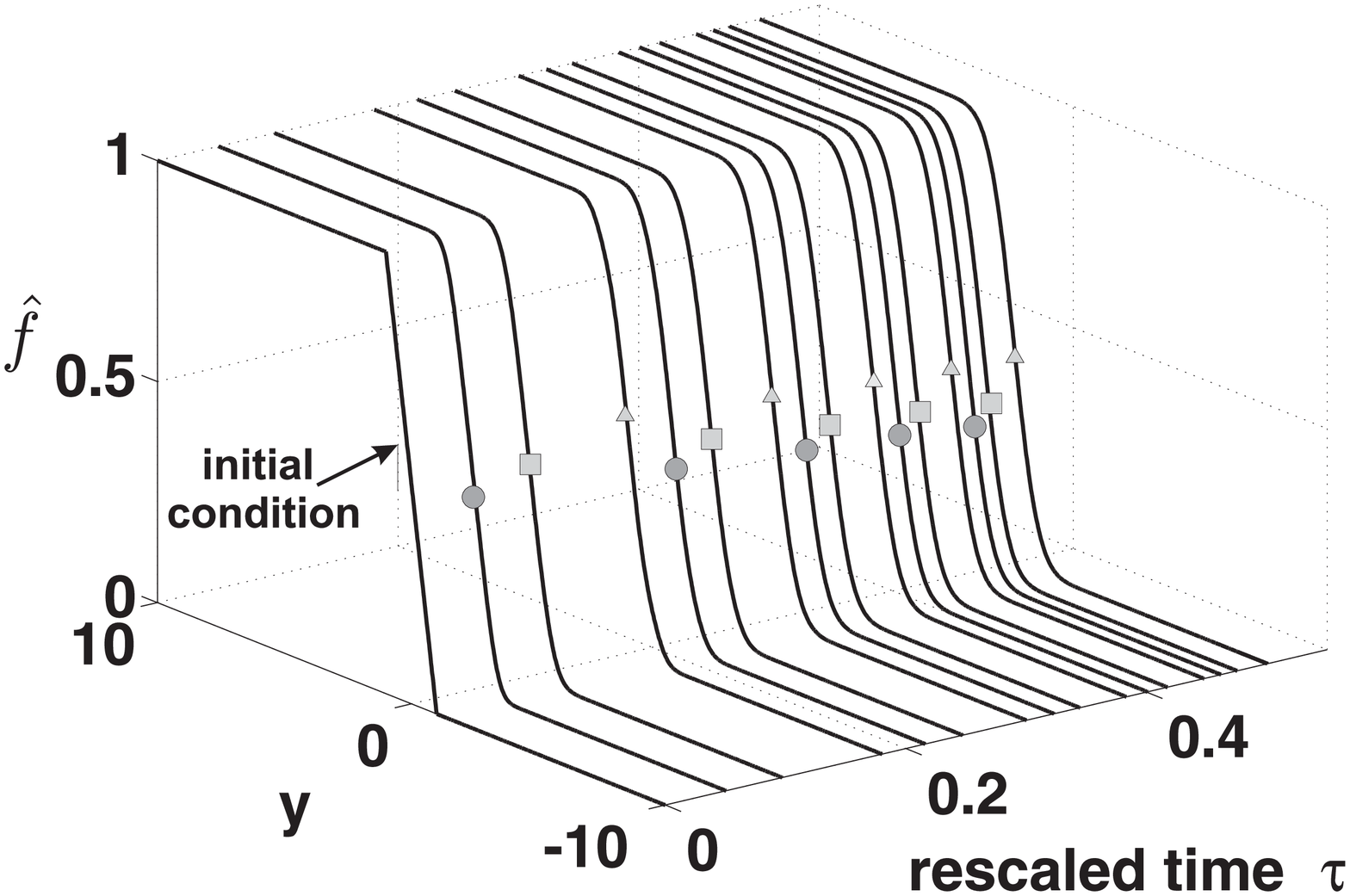} &
\includegraphics[width=0.49\linewidth]{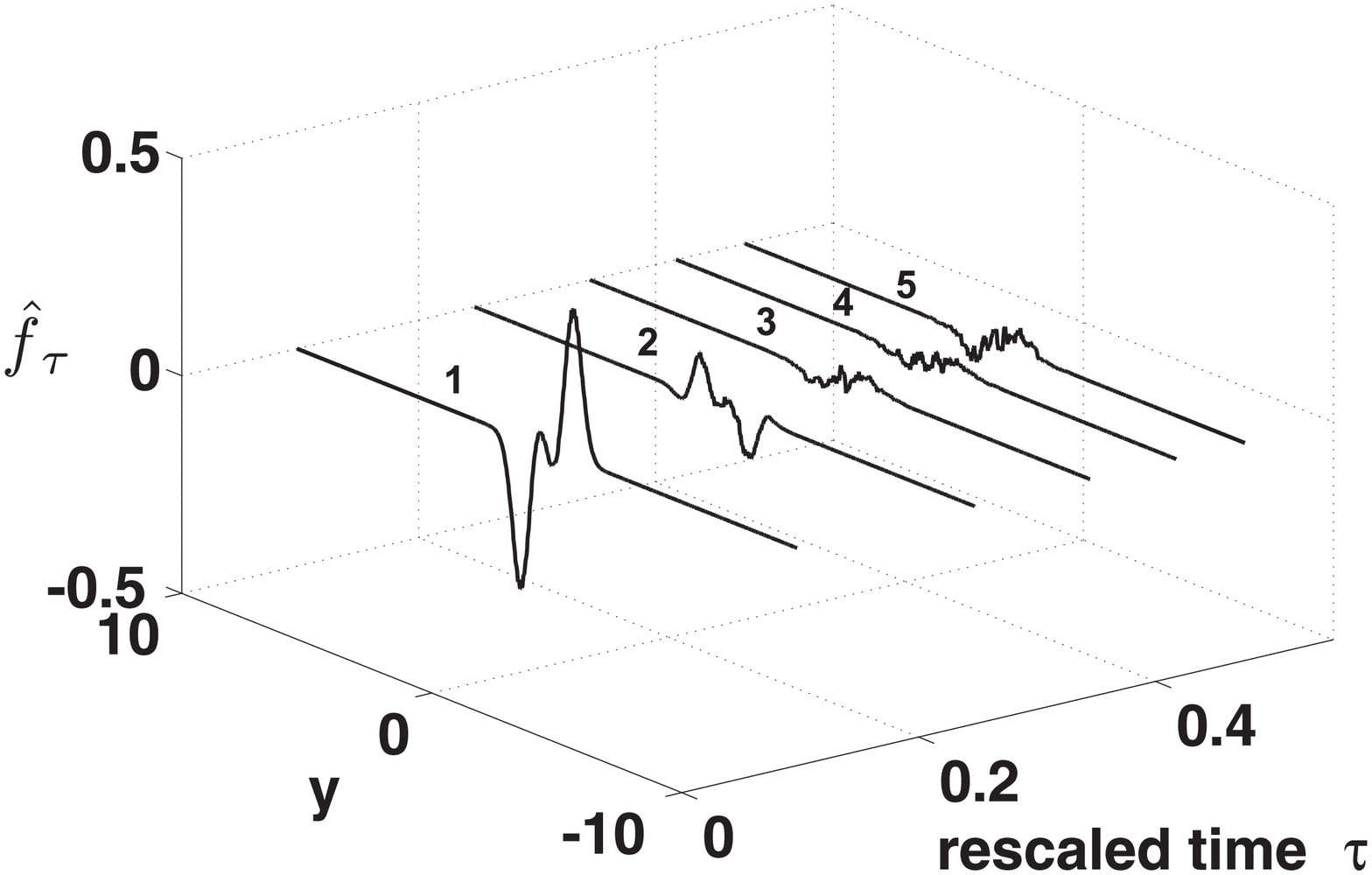} \\
    (a) & (b)
  \end{tabular}
  \end{center}
 \caption{Random walker simulation of 
 one-dimensional diffusion: (a) Coarse projective integration applied to the CDF
$\hat{f}$ in the dynamically co-evolving frame. The lines marked with circles correspond
to $\hat{f}(\tau_{1})$ solutions, the square-marked lines correspond to 
$\hat{f}(\tau_{2})$ solutions and the lines marked with triangles correspond to 
$\hat{f}$ obtained from projection at $\tau_{project}$.
 (b) Time derivative of $\hat{f}$ evaluated at the $1^{st}$--$5^{th}$ projective steps
of the co-evolving projective integration algorithm.}
 \label{hatfevolution}
 \end{figure}

\begin{figure}
 \begin{center}
 \begin{tabular} {cc}
 \includegraphics[width=0.49\linewidth]{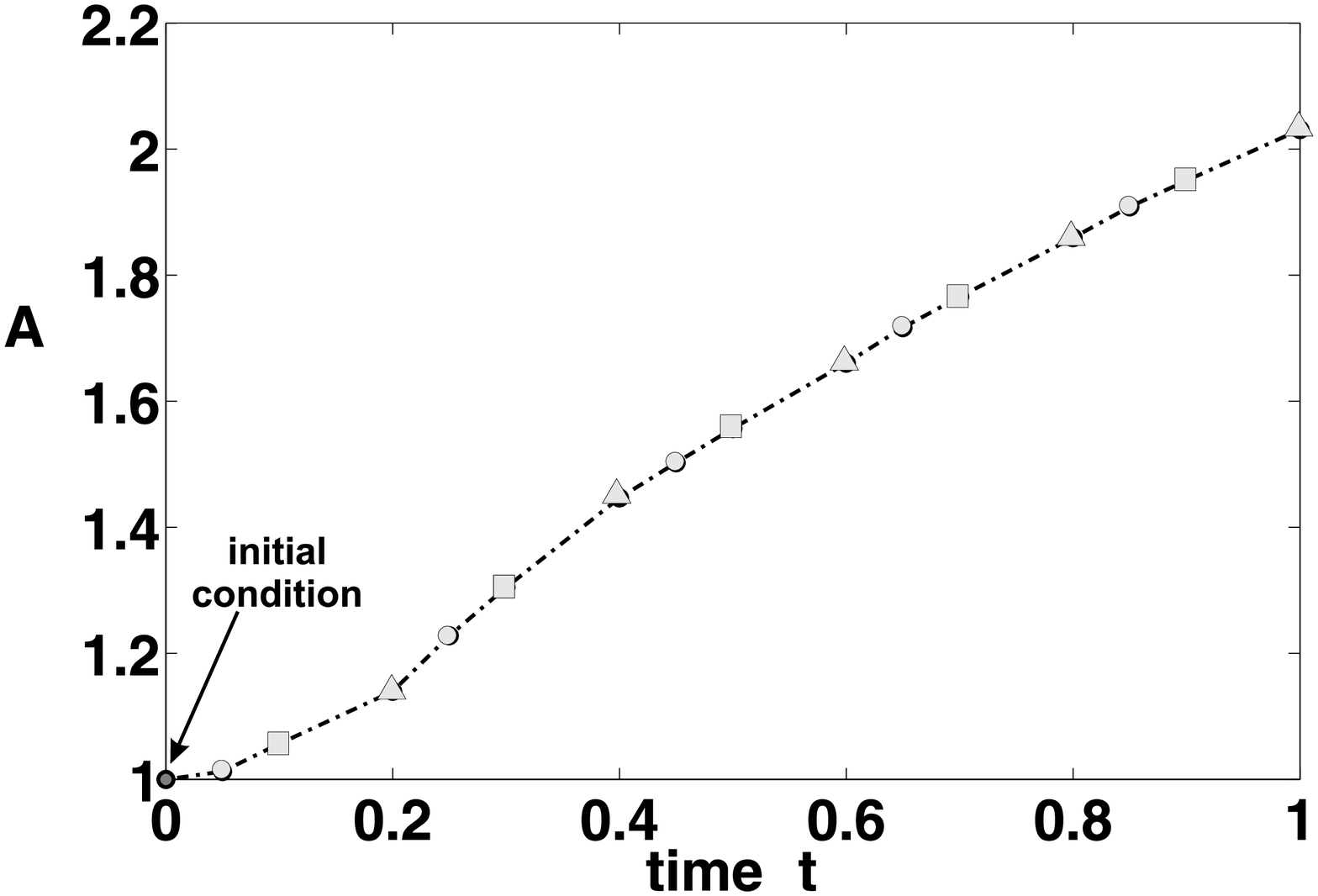} &
\includegraphics[width=0.49\linewidth]{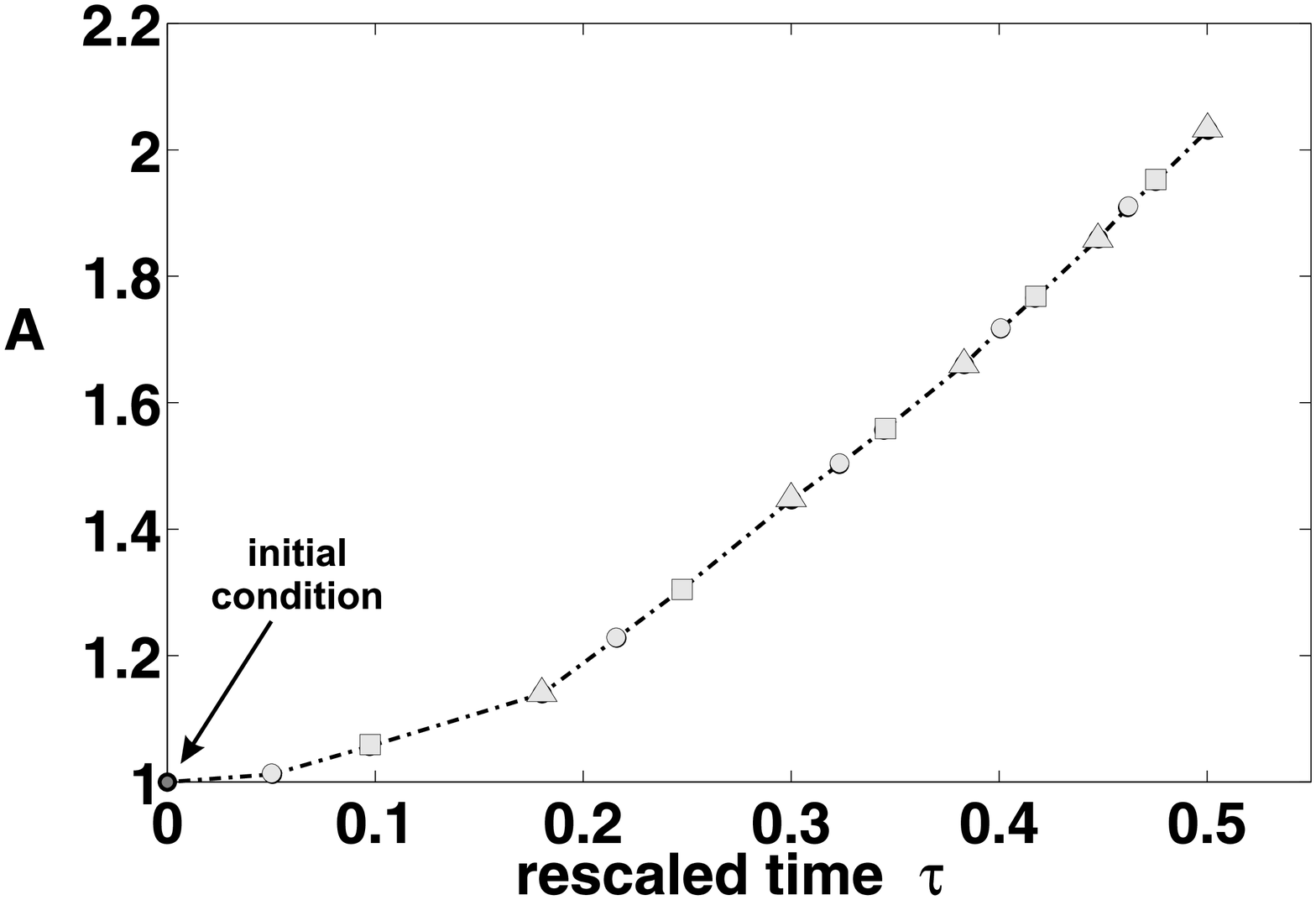} \\
    (a) & (b)
    \end{tabular}
  \end{center}
 \caption{Random walker simulation: Scale factor $A$ values computed from 
template condition (\ref{templateA}) during the application 
of the modified projective integration algorithm. The circle points correspond to
$A$-values computed at (a) $t_{1}$ -- reporting times (b) $\tau_{1}$ -- reporting times.
The square points correspond to (a) $t_{2}$ and (b) $\tau_{2}=\tau(t_{2})$ reporting times.
The $A$ values obtained at projection times $t_{project}$ ((b) $\tau_{project}$) are 
marked with triangles.} 
\label{scale_evolution}
 \end{figure}

The numerical results are shown in Figures \ref{hatfevolution} and \ref{scale_evolution}.
In Fig.\ref{hatfevolution} we present results of the first five 
coarse projective steps applied to the rescaled CDF $\hat{f}$.
The circle-marked lines correspond to reporting $\tau_{1}$-times, the square-marked lines to reporting 
$\tau_{2}$-times and the triangle-marked lines to projective $\tau_{project}$-times.
In Fig.\ref{scale_evolution} we plot the evolution of scale factor $A$
both in terms of time $t$ and of rescaled time $\tau$, as computed from the modified 
projective integration using the template condition (\ref{templateA}).

\section{A general approach for problems with asymptotic or approximate scale and
translational invariance}
\label{generalapproach}
In most real-world problems the issue of translation and rescaling
must be handled simultaneously. 
Even if the problem {\em is} self-similar,
our choice of templates for normalizing may lead to an
apparent translation with time. 
For example, the viscous Burgers equation
\begin{equation}
u_t = uu_x + \kappa u_{xx}
\end{equation}
has a self-similar solution
\begin{equation}
u(x,t) = t^{-1/2}w(xt^{-1/2})
\end{equation}
where $w(y)$ is not a symmetric function.  
Unfortunately if we
choose the wrong origin in the $y$ coordinate, we will find that
the evolving waveform is also traveling because we will actually
be looking at
\begin{equation}
u(x,t) = t^{-1/2}w((x-x_0)t^{-1/2} + x_0t^{-1/2})
\end{equation}
where the $x_0t^{-1/2}$ looks like a translation with time.

Many problems are only {\em asymptotically} self-similar.  
For example the solution of
\begin{equation}
 u_t = (1 + u^2)u_{xx}
 \end{equation}
asymptotically approaches the self-similar solution of the heat equation because
as $u$ decays the $u^2$ term becomes asymptotically small. 
Other problems may be {\em approximately} self-similar.  
For example, they may
approximately satisfy a scaling relationship such as
\begin{equation}
\label{cond1}
\operL_x\left(Bf\left(\frac{x}{A}\right)\right) = B^bA^a\operL_y\left(f\left(y\right)\right)[1 + {\rm O}(\varepsilon(t))]
\end{equation}
where $\varepsilon(t)$ is small and $y = x/A$.

Therefore, we want to look for time-dependent translations and
rescalings that yield a slowly varying waveform that can be
integrated in time more accurately because of its smaller time
derivatives.  
However, we cannot use the mechanism described in
the previous sections, where we examined ``purely'' scale invariant systems,
because we do not know $a$ and $b$ and they may not be
defined away from the asymptotic limit.

As the integration proceeds, it generates successive values of the
solution, $u(x,t)$. 
Periodically we need to apply three transformations to $u(x,t)$ to get
\begin{equation}
\label{3scale}
\hat{u}(y,t) = \frac{1}{B(t)}u(C(t) + A(t)y,t)
\end{equation}
with $y = [x - C(t)]/A(t)$
to {\em try} to get a $\hat{u}(y,t)$ which is {\it as independent of time $t$ as
possible}.  
Motivated by the approach for co-evolving systems that
are exactly self-similar, we will also consider a transformation of $t$ to
$\tau(t)$ in the expectation that $A(\tau)$ and $B(\tau)$
will evolve {\em approximately} exponentially in $\tau$ so that linear
projection of their logarithms will provide an accurate
approximation.

Template conditions are needed to determine the three scalings in
(\ref{3scale}).  
If we had a moderately good approximation to the
final waveform, it would be tempting to ask that the scalings be
chosen to minimize the difference between $\hat{u}(y,t)$ and that
approximate waveform.  
One might even choose the scalings so as to
minimize some norm of the time derivative.
However, such conditions are nonlinear, and
nonlinear conditions cause a problem in the projective step.  
A
projective step takes the form
\begin{equation}
\hat{u}_{project} = \hat{u}_2 + \frac{t_{project} -t_2}{t_2-t_1}[\hat{u}_2 - \hat{u}_1].
\end{equation}
In other words, $\hat{u}_{project}$ is a {\em linear} combination of $\hat{u}_2$ and
$\hat{u}_1$. 
If $\hat{u}_2$ and $\hat{u}_1$ satisfy a {\em linear} template
condition, $\hat{u}_{project}$ will also satisfy that condition automatically.  
That
is not necessarily true for a non-linear condition; after a
projection step one would have to re-apply the condition.

Therefore, we choose linear conditions to determine the scalings.
If the solution is non-negative - as it will be for many physically
based problems in which the variables we will use in the templates
are quantities such as density (represented at the microscopic level by
numbers of particles) there is a straightforward recipe.  
A shift can be determined by demanding that the center of gravity be
shifted to a fixed position, typically the origin.  
This is a
trivial calculation. 
For example, for uniform particles we simply
average their positions.  
In a continuum model we ask that
\begin{equation}
0 = {\int_{-\infty}^\infty y\hat{u}(y) \D y}
\end{equation}
which means that
\begin{equation}
\label{generalC}
C = \frac{\int_{-\infty}^\infty xu(x)\D x}{\int_{-\infty}^\infty u(x) \D x}.
\end{equation}
The $A$ scaling can be calculated by requiring that a certain fraction of the
integral of $\hat{u}$ lies within a specified interval, for example that
${\int_{-1}^1\hat{u}(y) \D y}=0.5{\int_{-\infty}^\infty\hat{u}(y) \D y}$,
although it is computationally easier to specify the second moment and require that
\begin{equation}
K = \frac{\int_{-\infty}^\infty y^2\hat{u}(y) \D y}{\int_{-\infty}^\infty\hat{u}(y) \D y}
\end{equation}
which implies that
\begin{equation}
\label{generalA}
A^2 = \frac{\int_{-\infty}^\infty x^2u(x+C)\D x}{K\int_{-\infty}^\infty u(x)\D x}.
\end{equation}
Note that the shift has been applied {\em before} the next moment is calculated.
The amplitude scaling $B$ can be calculated by requiring that the total mass,
\begin{equation}
\int_{-\infty}^\infty \hat{u}(y) \D y,
\end{equation}
remain constant, say equal to $\mu$, leading to
\begin{equation}
\label{generalB}
B = \frac{\int_{-\infty}^\infty u(x)\D x}{A\mu}.
\end{equation}
As the microscopic integration proceeds we calculate the values of $A$, $B$ and $C$ at
selected times $t_0$, $t_1$, $t_2,...$ and the rescaled solution
$\hat{u}(y,t_0), \hat{u}(y,t_1), \hat{u}(y,t_2), \cdots$ in the coarse variables.
We can decide whether it is appropriate to apply a projective step to the coarse solution based
on the local behavior of the coarse variables and the scaling values.

As we approach the region where the solution is close to self-similar, we need a way to compute
the transformation from $t$ to $\tau$ so that we can use exponential projective integration in $\tau$.
We can do this as follows.
We assume a form like (\ref{cond1}) in the PDE
\begin{equation}
\label{eq1}
u_t = \operL_x(u)
\end{equation}
and assume a solution of the form
\begin{equation}
\label{sol1}
u(x,t) = B(\tau)w\left(\frac{x}{A(\tau)},\tau\right)
\end{equation}
for some $\tau(t)$ and $w(y,\tau)$ which is slowly changing in $\tau$.
Then we get the equation
\begin{equation}
\label{selfeq}
\left(\frac{B_\tau}{B}w - \frac{A_\tau}{A}yw_y + w_\tau\right)\frac{\partial \tau}{\partial t} = B^{b-1}A^a\operL_yw[1+{\rm O}(\varepsilon(t))].
\end{equation}
We want to choose $A(\tau)$, $B(\tau)$, and $\tau(t)$ (which are completely at our choice) so that
{\em if there exists an approximately self-similar solution}, that is, a solution of
\begin{equation}
(b_1w - a_1yw_y)\frac{\partial \tau}{\partial t} = B^{b-1}A^a \operL_yw
\end{equation}
for some $a_1$ and $b_1$, $w$ tends to it.  
We
naturally choose $\tau_t = cB^{b-1}A^a$ for some constant $c$
and $A$ and $B$ so that $A_\tau/A$ and $B_\tau/B$ are nearly constant
and equal to $a_1$ and $b_1$ respectively.
(In practice, we will be choosing $A$ and $B$ to account for the observed growth in width and amplitudes
of the computed solution.)

If $A_\tau/A$ and $B_\tau/B$ are constant, then $A(\tau) = \exp(a_0 + a_1\tau)$ and
$B(\tau) = \exp(b_0 + b_1\tau)$.  Hence
\begin{equation}
\frac{\partial{t}}{\partial \tau} = \left({\frac{\partial \tau}{\partial t}}\right)^{-1} = \frac{\exp[-a(a_0 + a_1\tau)-(b-1)(b_0 + b_1\tau)]}{c}
\end{equation}
or
\begin{equation}
t = t_c -\frac{\exp[-a_0a-b_0(b-1)]}{(a_1a+b_1(b-1))c}\exp[-(a_1a+b_1(b-1))\tau].
\end{equation}
Note that we have this exponential behavior for $t(\tau)$ regardless of the actual values of $a$ and $b$.
Since the scale (and origin) of $\tau$ are arbitrary (we are picking them), we can rewrite this equation
as
\begin{equation}
\label{taudef}
t = t_c  + \beta\exp(\tau).
\end{equation}

Suppose now that we perform a calculation starting at $t_0$ and integrate to $t_1$ and $t_2$, 
computing $A(t_i)$ and
$B(t_i)$ as we proceed using template conditions.  
We can assume that $\tau_0 = 0$ since the origin is arbitrary.
We need to find $\tau_1$ and $\tau_2$.  
We have from (\ref{taudef})
\begin{equation}
(t_i - t_0)/\beta = \exp(\tau_i) - \exp(\tau_0)
\end{equation}
or, using $\tau_0 = 0$
\begin{equation}
\label{taurel}
\tau_i = \log(1 + (t_i-t_0)/\beta).
\end{equation}

Under the assumption that $A_\tau/A$ is constant we have
\begin{equation}
\label{arel}
\frac{1}{(\tau_1 - \tau_0)}\log\left(\frac{A_1}{A_0}\right) =\frac{1}{(\tau_2 - \tau_0)} \log\left(\frac{A_2}{A_0}\right)
\end{equation}
(and a similar relation for $B$).  Using (\ref{arel}), (\ref{taurel}) and $t_0 =\tau_0 = 0$, we get
\begin{equation}
\label{equend}
\log\left(1+\frac{t_1}{\beta}\right)\log\left(\frac{A_2}{A_0}\right) = \log\left(1+\frac{t_2}{\beta}\right)\log\left(\frac{A_1}{A_0}\right).
\end{equation}
This can be solved for $\beta$ and then we can use (\ref{taurel}) to find the $\tau_i$.
In the projective step, a suitable representation of the rescaled solution, $w$ 
can be projected in $\tau$, and
$A$ and $B$ are projected exponentially in $\tau$.

\subsection{An asymptotically scale and translationally invariant PDE example.}
\label{burgersexample}
\begin{figure}
\begin{center}
\includegraphics[width=0.6\linewidth]{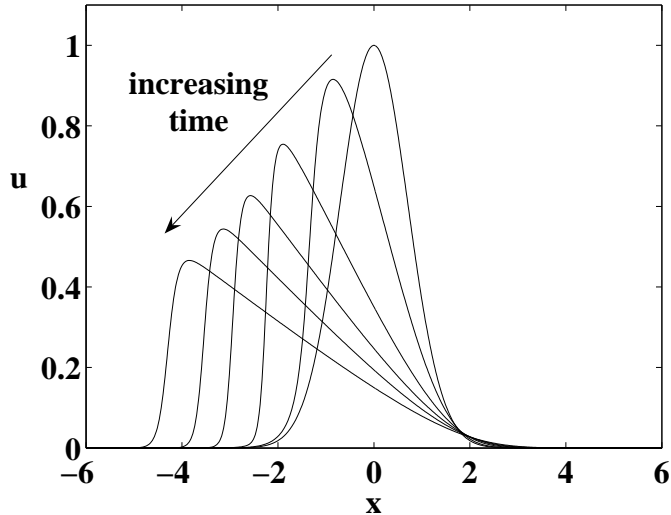}
\end{center}
\caption{Direct simulation of (\ref{Burgers}) with initial condition
$u(x,0)=\exp(-x^{2})$. The solution $u(x,t)$ 
travels both across scales and in space,
forming a steep interface which moves to the left part of the one-dimensional domain. Instances
of $u$ are presented at $t=0,1,3,5,7,10$.}
\label{Burgerevol}
\end{figure}

Below, we present some representative results from the application of this general
approach to the equation:
\begin{equation}
\label{Burgers}
u_{t}=\kappa (1+u^2) u_{xx} +u u_x \textrm{ with } \kappa =0.025
\end{equation}
which asymptotically approaches the solution of the viscous Burgers equation,
because as $u$ decreases the $u^2$ term becomes asymptotically small. 
Both translations
and rescalings have to be incorporated for this problem as one can see from the time
evolution of $u$ in Fig.\ref{Burgerevol}. 
Three template conditions are applied for the
evaluation of the shift $C$, and the scale factors $A$, $B$ at each reporting step, which have the form of
(\ref{generalC}), (\ref{generalA}) and (\ref{generalB}) respectively. 
The constants $K$ and $\mu$
appearing there are computed from the initial condition, i.e.:
\begin{equation}
K = \frac{\int_{-\infty}^\infty x^2 u(x,0) \D x}{\int_{-\infty}^\infty u(x,0) \D x}
\end{equation}
and the initial mass:
\begin{equation}
\mu=\int_{-\infty}^\infty u(x,0) \D x.
\end{equation}
The initial condition in our computations is $u(x,0)=\exp(-x^2)$. 
The applied boundary conditions are of Neummann type, the one-dimensional
computational domain $[-10,10]$ is discretized with $1001$ nodes and the spatial derivatives
of (\ref{Burgers}) are approximated with central finite differences. 

\begin{figure}
 \begin{center}
 \begin{tabular} {cc}
 \includegraphics[width=0.49\linewidth]{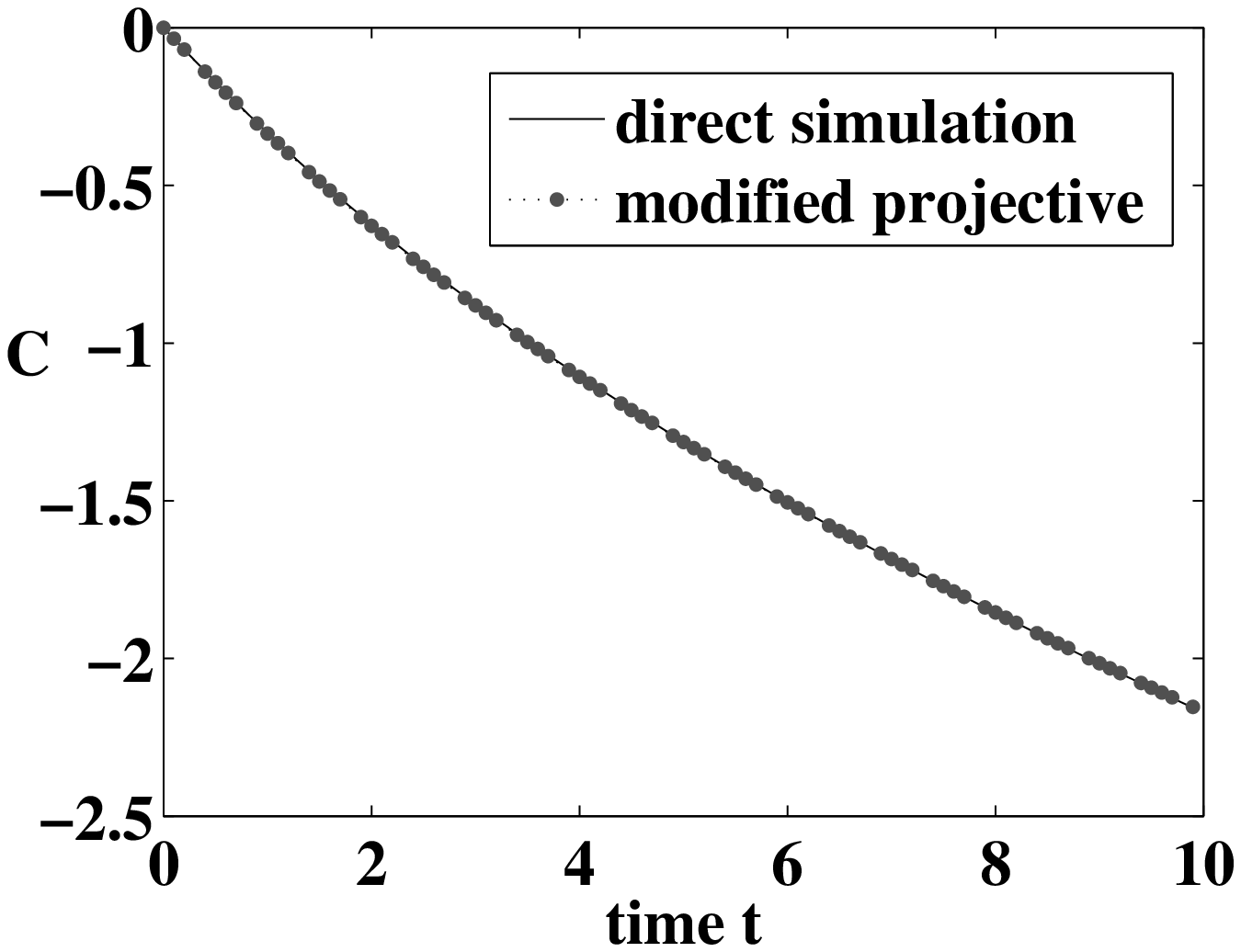} &
\includegraphics[width=0.49\linewidth]{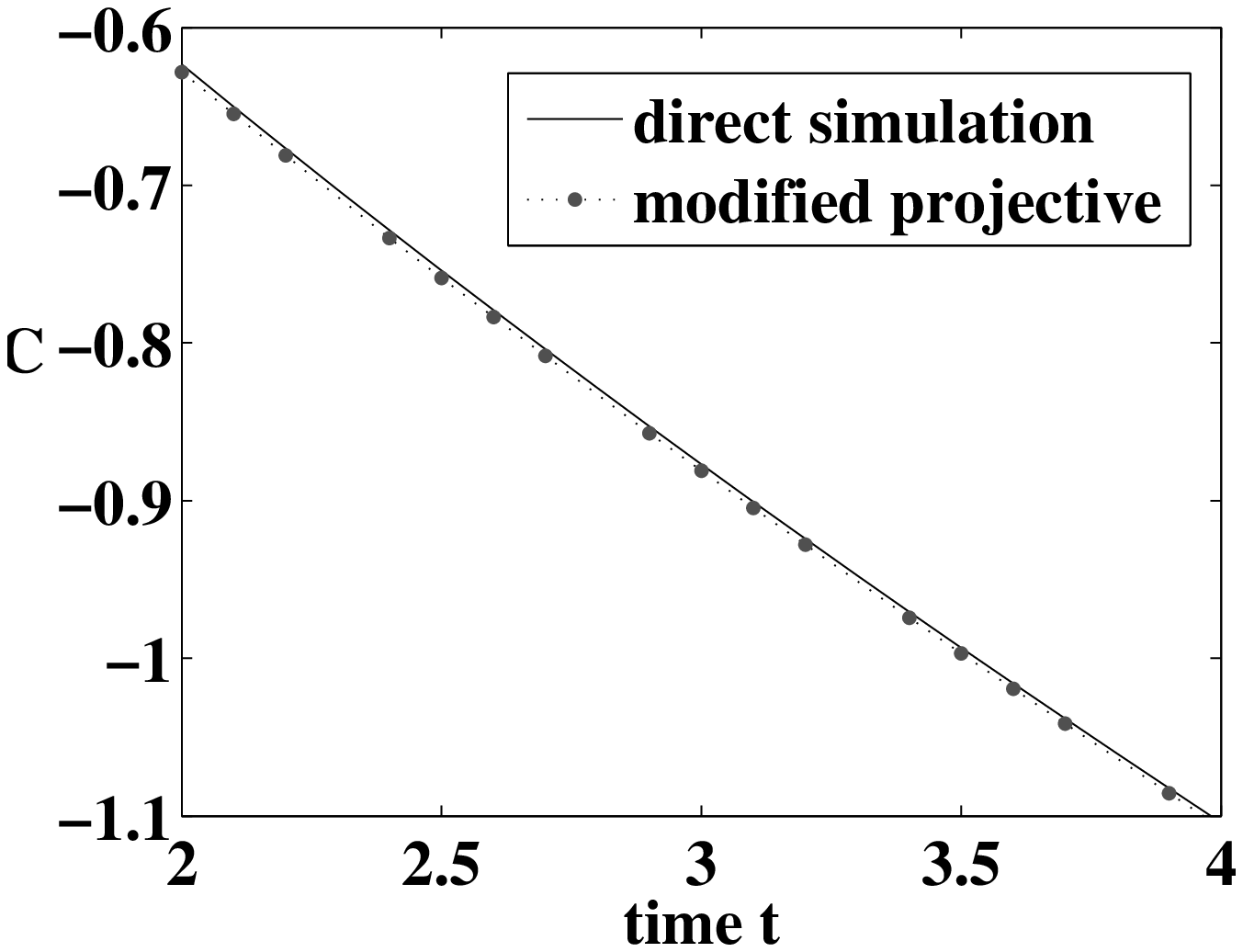}
\end{tabular}
\begin{tabular} {c}
(a)
\end{tabular}
 \begin{tabular} {cc}
 \includegraphics[width=0.49\linewidth]{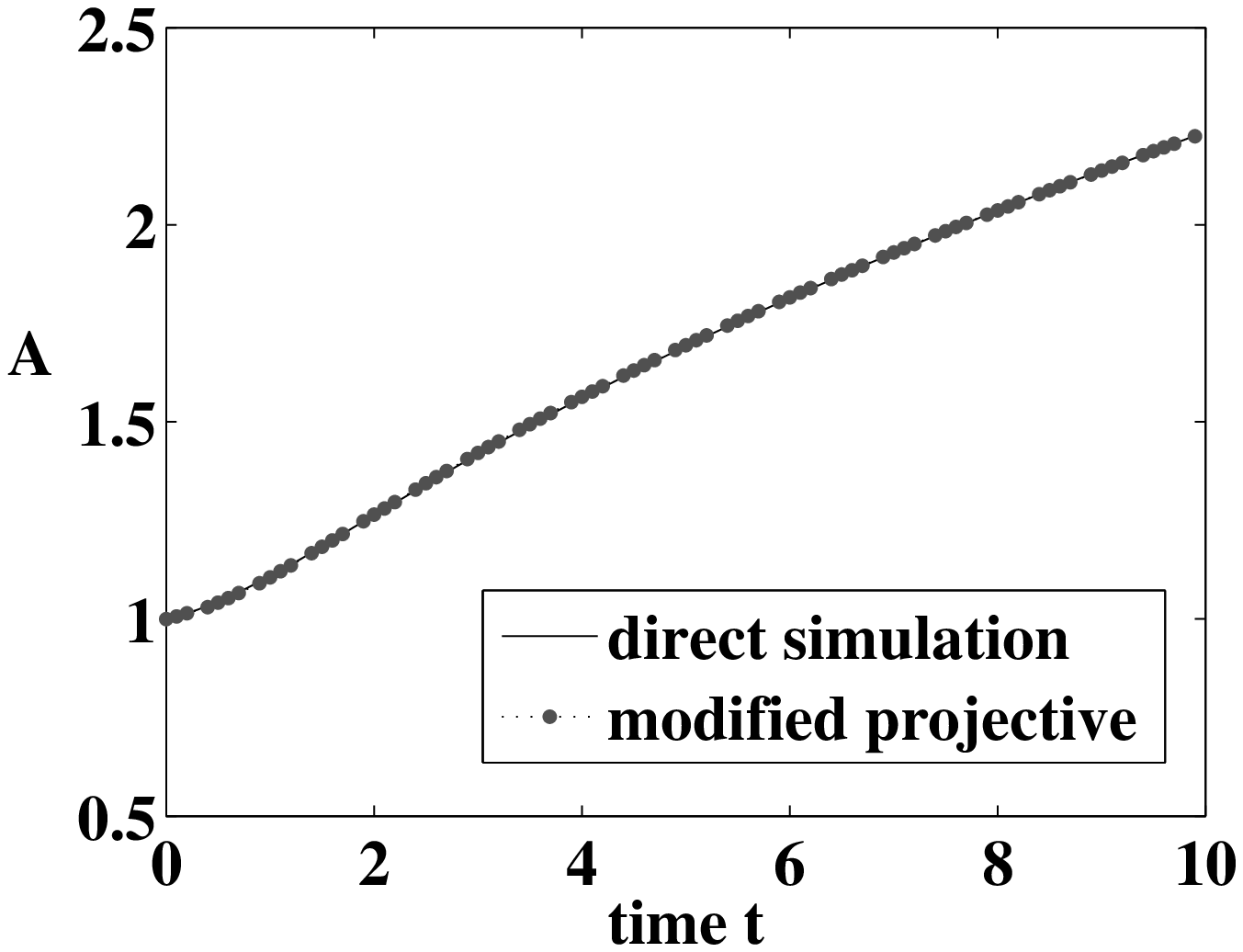} &
\includegraphics[width=0.49\linewidth]{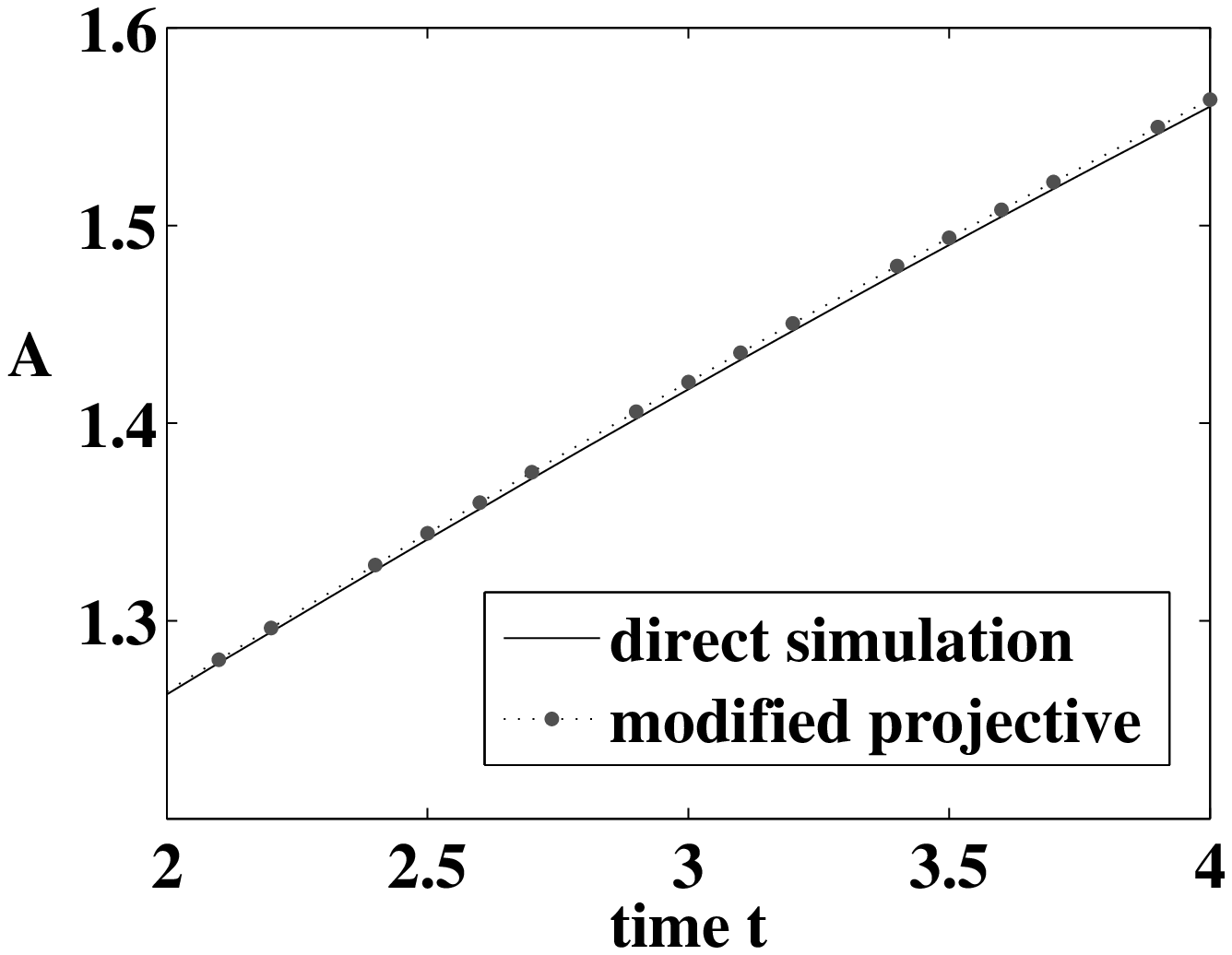}
\end{tabular}
\begin{tabular} {c}
(b)
\end{tabular}
 \begin{tabular} {cc}
 \includegraphics[width=0.49\linewidth]{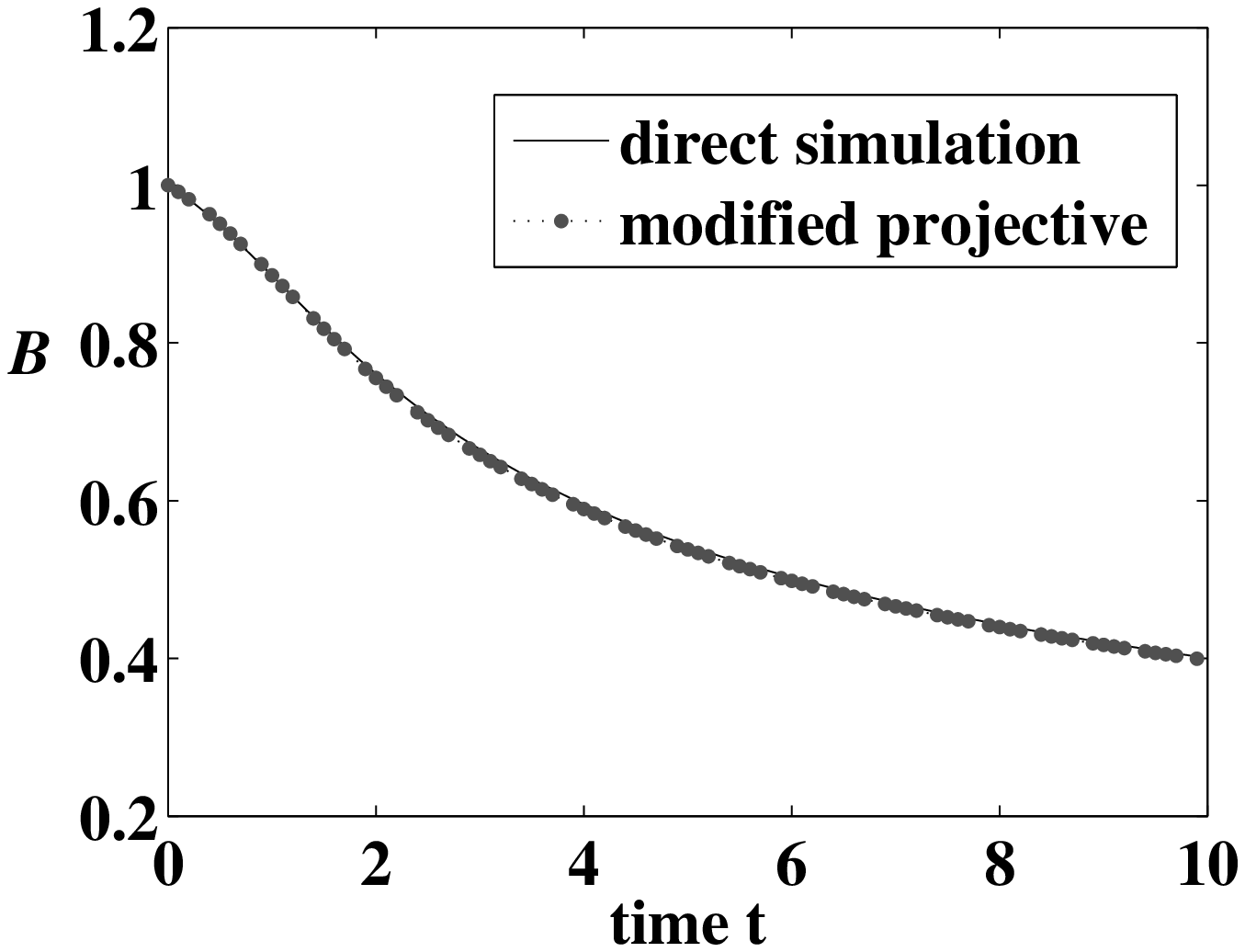} &
\includegraphics[width=0.49\linewidth]{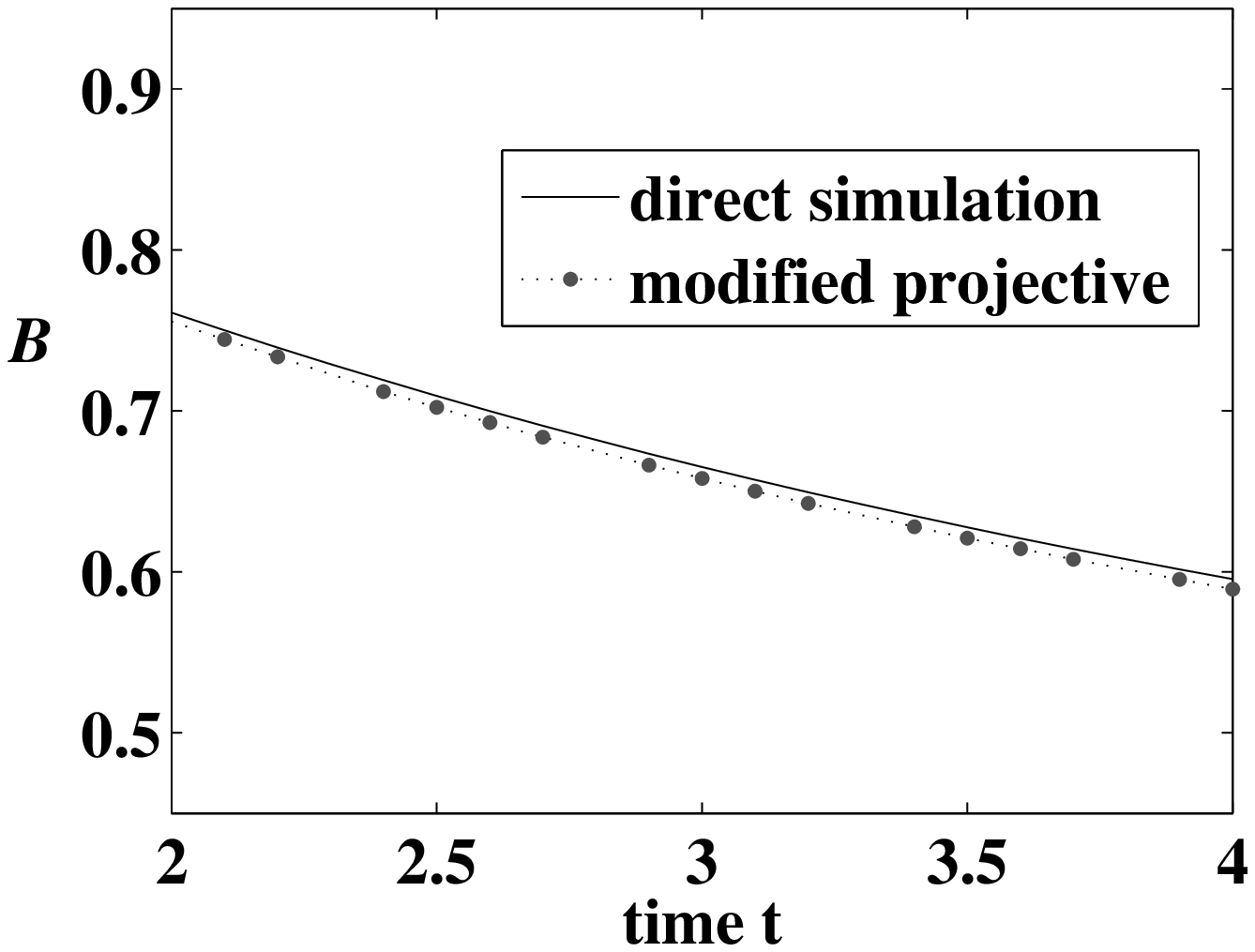}
\end{tabular}
\begin{tabular} {c}
(c)
\end{tabular}
  \end{center}
 \caption{ Computed values of (a) shift factor $C$ and scale factors (b) $A$, (c) $B$ 
from the application of template conditions (\ref{generalC}), 
(\ref{generalA}) and (\ref{generalB}) respectively, to the solution 
$u(x,t)$ of (\ref{Burgers}). The solid lines correspond to results
derived from direct simulation 
of (\ref{Burgers}). The $C$, $A$ and $B$ values, computed from the 
3-step template based projective method (see Section \ref{burgersexample}),
are depicted by dotted lines.}
 \label{ABCBurgers}
 \end{figure}

\begin{figure}
\begin{center}
\includegraphics[width=0.6\linewidth]{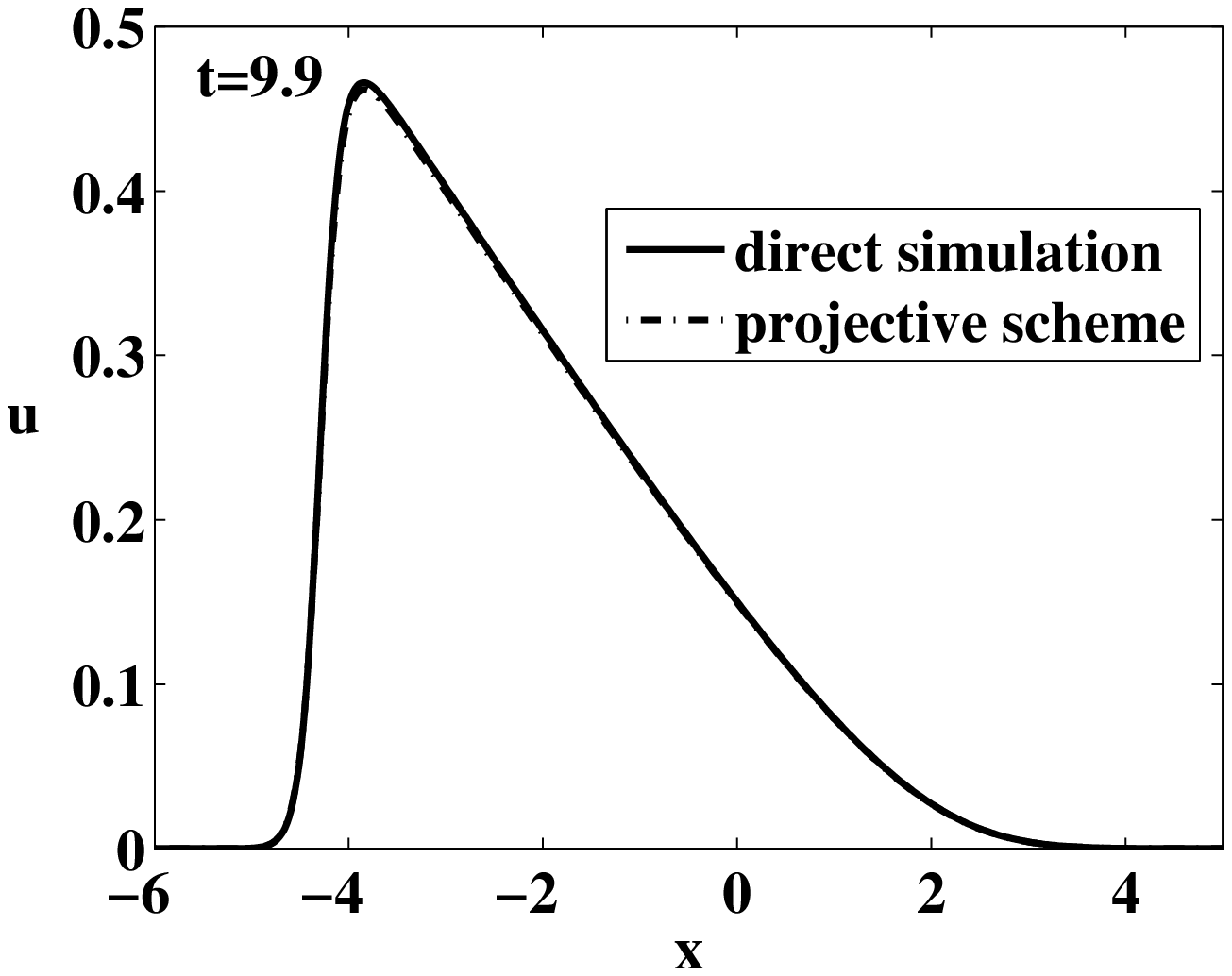}
\end{center}
\caption{Evolved solution $u$ of (\ref{Burgers}) at time $t=9.9$. The solid line represents
the solution derived from direct simulation and the dashed line corresponds to 
the solution computed by the 3-step template-based projective scheme described in Section
\ref{burgersexample}.}
\label{Burgersfinal}
\end{figure}

The direct time integration is
performed with an explicit Euler scheme, with $dt=10^{-5}$ ensuring the stability of the integration
for the given discretization.
The three reporting times are chosen so as $\Delta t =0.1$, while the projective step is taken
to be $T=0.2$.
At the initial stages of the integration, the solution is far from its asymptotically
self-similar shape solution and it is trivial to show that (\ref{Burgers}) is not scale invariant,
at least for large enough values of $u$. 
In this case we can project linearly in time
the shift factor $C$ the scale factors $A,B$ and  the renormalized solution $\hat{u}$
which is derived from the rescaling equation $\hat{u}(y)=1/B u(Ay+C)$. 
When the solution
approaches the self-similar regime, then we can apply the transformation of $t$ to $\tau$
and follow the procedure described in Sec. \ref{generalapproach}. 
The evolution of the factors
$A$ and $B$ computed from direct simulations and from the 3-step, template-based projective scheme
are depicted in Fig.\ref{ABCBurgers}. 
Finally, we illustrate the accuracy of this method,
presenting the solution computed from the direct simulation at time $t=9.9$ and the one
derived from the projective method at the same time (see Fig.\ref{Burgersfinal}).


\section{Summary and Conclusions}
\label{secdiscussion}

In this paper we have illustrated projective and coarse
projective integration
in a co-evolving frame for problems with continuous
symmetries, and in particular for problems with (coarse) scale 
invariance and
(coarse) translational invariance.
The system temporal evolution is observed in a traveling and/or 
dynamically renormalized frame, in which
transient solutions approaching traveling waves or self-similar 
solutions appear slowly changing,
and can be integrated more accurately because of
the smaller time derivatives.
Larger extrapolation time steps can thus
be applied without degrading the accuracy of the projected solution.
The simplest projective algorithm (projective forward Euler) was 
illustrated;
more sophisticated multistep and even implicit projective algorithms
are also possible (e.g. 
\cite{Gear:2002:CIB,Martinez:2004:CPK,Lee:2005:SOPI}).
We have illustrated several representative examples of template 
(pinning) conditions
that are used to dynamically define the coevolving frame in which 
projective integration
takes place.
Our model examples included both continuum and microscopic-based 
implementations.
The translationally invariant, co-traveling case was illustrated through 
the Nagumo
reaction-diffusion PDE \cite{Miura:1982:ACF,Murray:2002:MB,Beyn:2004:FSE}
in one spatial dimension, as well as an SSA-based 
\cite{Gillespie:1977:ESS} stochastic
implementation of the Nagumo kinetics for coarse projective integration.
The scale invariant case was illustrated through simple one-dimensional 
diffusion:
projective integration of the PDE version and coarse projective integration
of a stochastic implementation involving a large ensemble of random 
walkers were presented
and the results compared with direct, full simulation.
Finally, we described a more general projective method designed for 
systems with
{\it asymptotic} or even {\it approximate} invariance, and where scale 
invariance
and translational invariance co-exist.

The thrust of the paper was in describing and illustrating the methods, 
providing
some evidence and qualitative justification for the resulting computational
savings.
This constitutes only the starting point for the numerical analysis of 
the algorithms,
both for the deterministic and for the stochastic cases, which is the 
subject of
further research.
It is worth reiterating that the template based approach transforms 
traveling or self-similar
problems into steady-state ones; traveling wave speeds,
``scale velocities" $\xi_A, \xi_B$ and similarity exponents are simple 
and natural byproducts
of the approach.
Accelerating the computation of coarse self-similar shapes and coarse 
similarity exponents
for microscopic/stochastic simulators can be useful in a variety of 
disciplines, ranging
from microhydrodynamics to core collapse in star clusters \cite{Szell:2005:CCC}.
In the stochastic case, the accurate and efficient estimation of time 
derivatives from stochastic simulations
becomes a vital component of the algorithm, and one must move beyond simple
differencing and least squares estimators (like the ones we used here) 
to maximum likelihood
ones (see e.g. \cite{AitSahalia:2002:MLE}).
It is worth noting that -whether in the deterministic or in the 
coarse-grained case-
it is important to explore the relation between modern adaptive mesh 
techniques used for the
computation of self-similar solutions \cite{Budd:1999:SSN,Budd:1999:NSS,Ren:2000:IGR}
with the template-based ones presented in
\cite{Rowley:2003:RRS,Beyn:2004:FSE} and exploited here.
The most important factor in the success of coarse-grained projective 
integration (and of equation-free
computation in general) is the {\it lifting} step: the ability to 
construct ensembles of
fine-scale realizations consistent with a given macroscopic description.
For the problems discussed in this paper, the coarse observables in terms of
which the process was modeled allowed for a relatively easy
and computationally inexpensive lifting.
If the coarse-grained model is cast in terms of different coarse-grained 
observables
(e.g., particle pair correlation functions) the lifting step may become 
much more
difficult and expensive (see e.g. {\cite{Torquato:2002:RHM}).
Clearly, the cost of the lifting step (a very much problem-dependent 
feature) must be
factored in when evaluating the potential savings of coarse projective 
integration.
We close by reiterating that what we have presented is only a first step 
in the
study of coarse-grained projective integration algorithms for systems 
with (coarse)
continuous symmetries; the potential benefits illustrated here for two 
simple model problems
argue that the algorithms both on the continuum front (e.g. projective
integration for differential-algebraic equations \cite{Gear:2005:TEM}, 
relations to adaptive mesh algorithms)
and on the stochastic front (issues of lifting and estimation) warrant 
extensive further study.
\\

{\bf{Acknowledgements.}}  This work was partially supported by the Federal 
Fellowship Foundation of Greece and the NTUA through the Basic
Research Program ``Protagoras'' (MEK, AGB), by the Biotechnology and Biological
Sciences Research Council
and Linacre College, University of Oxford (RE), by the U.S. Department of Energy,
DARPA and a Guggenheim Fellowship (IGK).

\end{document}